\documentstyle[amstex,12pt, righttag,epsbox]{article}
\numberwithin{equation}{section}
\newcommand{\bysame}%
                   {\leavevmode\hbox to 3em{\hrulefill}\,}
\newtheorem{thm}{Theorem}[section]
\newtheorem{lem}{Lemma}[section]
\newtheorem{cor}{Corollary}[section]
\newtheorem{prop}{Proposition}[section]
\newtheorem{conj}{Conjecture}[section]
\newtheorem{defn}{Definition}[section]

\newtheorem{rem}{Remark}[section]

\newcommand{\qed}{q.e.d.}
\newcommand{\Fvertex}{F\left(\raisebox{-0.12cm}{
\begin{picture}(10,10)
\put(0,0){\psbox[scale=0.06]{V-insertionarrow}}
\end{picture}}\right)}
\newcommand{\Fcross}{F\left(\raisebox{-0.12cm}{
\begin{picture}(10,10)
\put(0,0){\psbox[scale=0.06]{crossarrow}}
\end{picture}}\right)}
\newcommand{\Fccross}{F\left(\raisebox{-0.12cm}{
\begin{picture}(10,10)
\put(0,0){\psbox[scale=0.06]{C-insertionarrow}}
\end{picture}}\right)}
\newcommand{\Ftvertex}{F^t\left(\raisebox{-0.12cm}{
\begin{picture}(10,10)
\put(0,0){\psbox[scale=0.06]{V-insertionarrow}}
\end{picture}}\right)}
\newcommand{\Fpcross}{F\left(\raisebox{-0.12cm}{
\begin{picture}(10,10)
\put(0,0){\psbox[scale=0.06]{pcrossarrow}}
\end{picture}}\right)}
\newcommand{\Fncross}{F\left(\raisebox{-0.12cm}{
\begin{picture}(10,10)
\put(0,0){\psbox[scale=0.06]{ncrossarrow}}
\end{picture}}\right)}
\newcommand{\Lcircleb}{\left(L\raisebox{-0.18cm}{
\begin{picture}(12,11)
\put(0,0){\psbox[scale=0.06]{Circle}}
\end{picture}}\right)}
\newcommand{\Lpgraph}{\left(L\raisebox{-0.2cm}{
\begin{picture}(12,11)
\put(0,0){\psbox[scale=0.06]{graph-a4}}
\end{picture}}\right)}
\newcommand{\Lngraph}{\left(L\raisebox{-0.2cm}{
\begin{picture}(12,11)
\put(0,0){\psbox[scale=0.06]{graph-a5}}
\end{picture}}\right)}
\newcommand{\Lgraph}{\left(L\raisebox{-0.16cm}{
\begin{picture}(16,10)
\put(0,0){\psbox[scale=0.06]{graph-a3}}
\end{picture}}\right)}
\newcommand{\LLjcrossb}{\left(L^{(j,0)}\raisebox{-0.17cm}{
\begin{picture}(10,10)
\put(0,0){\psbox[scale=0.06]{crossarrow}}
\end{picture}}\right)}
\newcommand{\Lcroossb}{\left(L\raisebox{-0.17cm}{
\begin{picture}(10,10)
\put(0,0){\psbox[scale=0.06]{croossarrow}}
\end{picture}}\right)}
\newcommand{\Ljcroossb}{\left(L^{(j)}\raisebox{-0.17cm}{
\begin{picture}(10,10)
\put(0,0){\psbox[scale=0.06]{croossarrow}}
\end{picture}}\right)}
\newcommand{\Lonecroossb}{\left(L^{(1)}\raisebox{-0.17cm}{
\begin{picture}(10,10)
\put(0,0){\psbox[scale=0.06]{croossarrow}}
\end{picture}}\right)}
\newcommand{\Pcrossb}{\left(\raisebox{-0.22cm}{
\begin{picture}(16,16)
\put(0,0){\psbox[scale=0.1]{pcrossarrow}}
\end{picture}}\right)}
\newcommand{\Lpcrossb}{\left(L\raisebox{-0.17cm}{
\begin{picture}(10,10)
\put(0,0){\psbox[scale=0.06]{pcrossarrow}}
\end{picture}}\right)}
\newcommand{\Lzpcrossb}{\left(L^{(0)}\raisebox{-0.17cm}{
\begin{picture}(10,10)
\put(0,0){\psbox[scale=0.06]{pcrossarrow}}
\end{picture}}\right)}
\newcommand{\LLjpcrossb}{\left(L^{(j-1,0)}\raisebox{-0.17cm}{
\begin{picture}(10,10)
\put(0,0){\psbox[scale=0.06]{pcrossarrow}}
\end{picture}}\right)}
\newcommand{\Lncrossb}{\left(L\raisebox{-0.17cm}{
\begin{picture}(10,10)
\put(0,0){\psbox[scale=0.06]{ncrossarrow}}
\end{picture}}\right)}
\newcommand{\Ncrossb}{\left(\raisebox{-0.22cm}{
\begin{picture}(16,16)
\put(0,0){\psbox[scale=0.1]{ncrossarrow}}
\end{picture}}\right)}
\newcommand{\Lzncrossb}{\left(L^{(0)}\raisebox{-0.17cm}{
\begin{picture}(10,10)
\put(0,0){\psbox[scale=0.06]{ncrossarrow}}
\end{picture}}\right)}
\newcommand{\LLjncrossb}{\left(L^{(j-1,0)}\raisebox{-0.17cm}{
\begin{picture}(10,10)
\put(0,0){\psbox[scale=0.06]{ncrossarrow}}
\end{picture}}\right)}
\newcommand{\Unfoldingb}{\left(\raisebox{-0.22cm}{
\begin{picture}(16,16)
\put(0,0){\psbox[scale=0.1]{unfoldingarrow}}
\end{picture}}\right)}
\newcommand{\Lunfoldingb}{\left(L\raisebox{-0.17cm}{
\begin{picture}(10,10)
\put(0,0){\psbox[scale=0.06]{unfoldingarrow}}
\end{picture}}\right)}
\newcommand{\Unfoldingsssb}{\left(\raisebox{-0.22cm}{
\begin{picture}(16,16)
\put(0,0){\psbox[scale=0.1]{unfoldingarrow4}}
\end{picture}}\right)}
\newcommand{\Unfoldb}{\left(\raisebox{-0.22cm}{
\begin{picture}(16,16)
\put(0,0){\psbox[scale=0.1]{unfoldingarrow5}}
\end{picture}}\right)}
\newcommand{\Lonecasimirb}{\left(L^{(1)}\raisebox{-0.17cm}{
\begin{picture}(10,10)
\put(0,0){\psbox[scale=0.06]{C-insertionarrow}}
\end{picture}}\right)}
\newcommand{\Pcrosssb}{\left(\raisebox{-0.22cm}{
\begin{picture}(16,16)
\put(0,0){\psbox[scale=0.1]{pcrossarrow2}}
\end{picture}}\right)}
\newcommand{\Dpcrossb}{\left(\raisebox{-0.73cm}{
\begin{picture}(35,44)
\put(0,0){\psbox[scale=0.17]{D-pcross}}
\end{picture}}\right)}
\newcommand{\Dncrossb}{\left(\raisebox{-0.73cm}{
\begin{picture}(35,44)
\put(0,0){\psbox[scale=0.17]{D-ncross}}
\end{picture}}\right)}
\newcommand{\Dunfoldingb}{\left(\raisebox{-0.73cm}{
\begin{picture}(34,44)
\put(0,0){\psbox[scale=0.17]{D-unfolding}}
\end{picture}}\right)}
\newcommand{\Ncrosssb}{\left(\raisebox{-0.22cm}{
\begin{picture}(16,16)
\put(0,0){\psbox[scale=0.1]{ncrossarrow2}}
\end{picture}}\right)}
\newcommand{\DDncrossb}{\left(\raisebox{-0.56cm}{
\begin{picture}(44,34)
\put(0,0){\psbox[scale=0.17]{DD-ncross}}
\end{picture}}\right)}
\newcommand{\DDpcrossb}{\left(\raisebox{-0.56cm}{
\begin{picture}(44,34)
\put(0,0){\psbox[scale=0.17]{DD-pcross}}
\end{picture}}\right)}
\newcommand{\Crosssb}{\left(\raisebox{-0.63cm}{
\begin{picture}(25,37)
\put(0,0){\psbox[scale=0.14]{crossarrow2}}
\end{picture}}\right)}
\newcommand{\Crossssb}{\left(\raisebox{-0.61cm}{
\begin{picture}(25,36)
\put(0,0){\psbox[scale=0.14]{crossarrow3}}
\end{picture}}\right)}
\newcommand{\Unfoldinb}{\left(\raisebox{-0.64cm}{
\begin{picture}(24,38)
\put(0,0){\psbox[scale=0.14]{unfoldingarrow0}}
\end{picture}}\right)}
\newcommand{\Unfoldinnb}{\left(\raisebox{-0.61cm}{
\begin{picture}(24,36)
\put(0,0){\psbox[scale=0.14]{unfoldingarrow6}}
\end{picture}}\right)}
\newcommand{\Unfoldingsb}{\left(\raisebox{-0.64cm}{
\begin{picture}(28,38)
\put(0,0){\psbox[scale=0.14]{unfoldingarrow2}}
\end{picture}}\right)}
\newcommand{\Unfoldingssb}{\left(\raisebox{-0.66cm}{
\begin{picture}(28,39)
\put(0,0){\psbox[scale=0.14]{unfoldingarrow3}}
\end{picture}}\right)}
\newcommand{\Graphaa}{\left(\raisebox{-0.44cm}{
\begin{picture}(21,28)
\put(0,0){\psbox[scale=0.1]{graph-a1}}
\end{picture}}\right)}
\newcommand{\Graphab}{\left(\raisebox{-0.44cm}{
\begin{picture}(21,28)
\put(0,0){\psbox[scale=0.1]{graph-a2}}
\end{picture}}\right)}
\newcommand{\Graphac}{\left(\raisebox{-0.28cm}{
\begin{picture}(26,16)
\put(0,0){\psbox[scale=0.1]{graph-a3}}
\end{picture}}\right)}
\newcommand{\Graphap}{\left(\raisebox{-0.44cm}{
\begin{picture}(21,28)
\put(0,0){\psbox[scale=0.1]{graph-a4}}
\end{picture}}\right)}
\newcommand{\Graphan}{\left(\raisebox{-0.44cm}{
\begin{picture}(21,28)
\put(0,0){\psbox[scale=0.1]{graph-a5}}
\end{picture}}\right)}
\newcommand{\Graphba}{\left(\raisebox{-0.25cm}{
\begin{picture}(29,18)
\put(0,0){\psbox[scale=0.1]{graph-b1}}
\end{picture}}\right)}
\newcommand{\Graphbb}{\left(\raisebox{-0.25cm}{
\begin{picture}(29,18)
\put(0,0){\psbox[scale=0.1]{graph-b2}}
\end{picture}}\right)}
\newcommand{\Graphbc}{\left(\raisebox{-0.25cm}{
\begin{picture}(29,18)
\put(0,0){\psbox[scale=0.1]{graph-b3}}
\end{picture}}\right)}
\newcommand{\Graphbd}{\left(\raisebox{-0.25cm}{
\begin{picture}(29,18)
\put(0,0){\psbox[scale=0.1]{graph-b5}}
\end{picture}}\right)}
\newcommand{\Graphbe}{\left(\raisebox{-0.25cm}{
\begin{picture}(29,18)
\put(0,0){\psbox[scale=0.1]{graph-b6}}
\end{picture}}\right)}
\newcommand{\Circleb}{\left(\raisebox{-0.25cm}{
\begin{picture}(18.67,17.57)
\put(0,0){\psbox[scale=0.1]{Circle}}
\end{picture}}\right)}
\newcommand{\scross}{\raisebox{-0.17cm}{
\begin{picture}(10,10)
\put(0,0){\psbox[scale=0.06]{crossarrow}}
\end{picture}}}
\newcommand{\scrooss}{\raisebox{-0.17cm}{
\begin{picture}(10,10)
\put(0,0){\psbox[scale=0.06]{croossarrow}}
\end{picture}}}
\newcommand{\spcross}{\raisebox{-0.17cm}{
\begin{picture}(10,10)
\put(0,0){\psbox[scale=0.06]{pcrossarrow}}
\end{picture}}}
\newcommand{\sncross}{\raisebox{-0.17cm}{
\begin{picture}(10,10)
\put(0,0){\psbox[scale=0.06]{ncrossarrow}}
\end{picture}}}
\newcommand{\sccross}{\raisebox{-0.17cm}{
\begin{picture}(10,10)
\put(0,0){\psbox[scale=0.06]{C-insertionarrow}}
\end{picture}}}
\newcommand{\sunfold}{\raisebox{-0.17cm}{
\begin{picture}(10,10)
\put(0,0){\psbox[scale=0.06]{unfoldingarrow}}
\end{picture}}}
\newcommand{\suunfold}{\raisebox{-0.17cm}{
\begin{picture}(10,10)
\put(0,0){\psbox[scale=0.06]{unfoldingarrow4}}
\end{picture}}}
\newcommand{\stunfold}{\raisebox{-0.17cm}{
\begin{picture}(10,10)
\put(0,0){\psbox[scale=0.06]{unfoldingarrow5}}
\end{picture}}}
\newcommand{\South}{\raisebox{-0.6cm}{
\begin{picture}(41,36)
\put(0,0){\psbox[scale=0.17]{South}}
\end{picture}}}
\newcommand{\North}{\raisebox{-0.6cm}{
\begin{picture}(41,36)
\put(0,0){\psbox[scale=0.17]{North}}
\end{picture}}}
\newcommand{\East}{\raisebox{-0.6cm}{
\begin{picture}(41,36)
\put(0,0){\psbox[scale=0.17]{East}}
\end{picture}}}
\newcommand{\West}{\raisebox{-0.6cm}{
\begin{picture}(41,36)
\put(0,0){\psbox[scale=0.17]{West}}
\end{picture}}}
\newcommand{\triple}{\raisebox{-0.6cm}{
\begin{picture}(41,36)
\put(0,0){\psbox[scale=0.17]{triple}}
\end{picture}}}
\newcommand{\Southb}{\raisebox{-0.6cm}{
\begin{picture}(41,36)
\put(0,0){\psbox[scale=0.17]{Southb}}
\end{picture}}}
\newcommand{\Northb}{\raisebox{-0.6cm}{
\begin{picture}(41,36)
\put(0,0){\psbox[scale=0.17]{Northb}}
\end{picture}}}
\newcommand{\Eastb}{\raisebox{-0.6cm}{
\begin{picture}(41,36)
\put(0,0){\psbox[scale=0.17]{Eastb}}
\end{picture}}}
\newcommand{\Westb}{\raisebox{-0.6cm}{
\begin{picture}(41,36)
\put(0,0){\psbox[scale=0.17]{Westb}}
\end{picture}}}
\newcommand{\trib}{\raisebox{-0.6cm}{
\begin{picture}(41,36)
\put(0,0){\psbox[scale=0.17]{tri}}
\end{picture}}}
\newcommand{\tripleabc}{\raisebox{-0.6cm}{
\begin{picture}(41,36)
\put(0,0){\psbox[scale=0.17]{triple-abc}}
\end{picture}}}
\newcommand{\triplea}{\raisebox{-0.6cm}{
\begin{picture}(41,36)
\put(0,0){\psbox[scale=0.17]{triple-a}}
\end{picture}}}
\newcommand{\tripleb}{\raisebox{-0.6cm}{
\begin{picture}(41,36)
\put(0,0){\psbox[scale=0.17]{triple-b}}
\end{picture}}}
\newcommand{\tripled}{\raisebox{-0.6cm}{
\begin{picture}(41,36)
\put(0,0){\psbox[scale=0.17]{triple-d}}
\end{picture}}}
\newcommand{\triplee}{\raisebox{-0.6cm}{
\begin{picture}(41,36)
\put(0,0){\psbox[scale=0.17]{triple-e}}
\end{picture}}}
\newcommand{\triplef}{\raisebox{-0.6cm}{
\begin{picture}(41,36)
\put(0,0){\psbox[scale=0.17]{triple-f}}
\end{picture}}}
\newcommand{\tripleg}{\raisebox{-0.6cm}{
\begin{picture}(41,36)
\put(0,0){\psbox[scale=0.17]{triple-g}}
\end{picture}}}
\newcommand{\tripleh}{\raisebox{-0.6cm}{
\begin{picture}(41,36)
\put(0,0){\psbox[scale=0.17]{triple-h}}
\end{picture}}}
\newcommand{\triplen}{\raisebox{-0.6cm}{
\begin{picture}(41,36)
\put(0,0){\psbox[scale=0.17]{triple-n}}
\end{picture}}}
\newcommand{\evv}{\raisebox{-0.6cm}{
\begin{picture}(48,37)
\put(0,0){\psbox[scale=0.17]{EastVV}}
\end{picture}}}
\newcommand{\evvs}{\raisebox{-0.6cm}{
\begin{picture}(48,37)
\put(0,0){\psbox[scale=0.17]{EastVVS}}
\end{picture}}}
\newcommand{\evvc}{\raisebox{-0.6cm}{
\begin{picture}(48,37)
\put(0,0){\psbox[scale=0.17]{EastVVC}}
\end{picture}}}
\newcommand{\wvv}{\raisebox{-0.6cm}{
\begin{picture}(48,37)
\put(0,0){\psbox[scale=0.17]{WestVV}}
\end{picture}}}
\newcommand{\wvvc}{\raisebox{-0.6cm}{
\begin{picture}(48,37)
\put(0,0){\psbox[scale=0.17]{WestVVC}}
\end{picture}}}
\newcommand{\wvvs}{\raisebox{-0.6cm}{
\begin{picture}(48,37)
\put(0,0){\psbox[scale=0.17]{WestVVS}}
\end{picture}}}
\newcommand{\svv}{\raisebox{-0.6cm}{
\begin{picture}(48,37)
\put(0,0){\psbox[scale=0.17]{SouthVV}}
\end{picture}}}
\newcommand{\svvs}{\raisebox{-0.6cm}{
\begin{picture}(48,37)
\put(0,0){\psbox[scale=0.17]{SouthVVS}}
\end{picture}}}
\newcommand{\svvc}{\raisebox{-0.6cm}{
\begin{picture}(48,37)
\put(0,0){\psbox[scale=0.17]{SouthVVC}}
\end{picture}}}
\newcommand{\nvv}{\raisebox{-0.6cm}{
\begin{picture}(48,37)
\put(0,0){\psbox[scale=0.17]{NorthVV}}
\end{picture}}}
\newcommand{\nvvs}{\raisebox{-0.6cm}{
\begin{picture}(48,37)
\put(0,0){\psbox[scale=0.17]{NorthVVS}}
\end{picture}}}
\newcommand{\nvvc}{\raisebox{-0.6cm}{
\begin{picture}(48,37)
\put(0,0){\psbox[scale=0.17]{NorthVVC}}
\end{picture}}}
\newcommand{\nvvb}{\raisebox{-0.6cm}{
\begin{picture}(48,37)
\put(0,0){\psbox[scale=0.17]{NVV}}
\end{picture}}}
\newcommand{\nvvsb}{\raisebox{-0.6cm}{
\begin{picture}(41,36)
\put(0,0){\psbox[scale=0.17]{NVVS2}}
\end{picture}}}
\newcommand{\nvvcb}{\raisebox{-0.6cm}{
\begin{picture}(41,36)
\put(0,0){\psbox[scale=0.17]{NVVC2}}
\end{picture}}}
\newcommand{\svvb}{\raisebox{-0.6cm}{
\begin{picture}(48,37)
\put(0,0){\psbox[scale=0.17]{SVV}}
\end{picture}}}
\newcommand{\svvsb}{\raisebox{-0.6cm}{
\begin{picture}(48,37)
\put(0,0){\psbox[scale=0.17]{SVVS2}}
\end{picture}}}
\newcommand{\svvcb}{\raisebox{-0.6cm}{
\begin{picture}(48,37)
\put(0,0){\psbox[scale=0.17]{SVVC2}}
\end{picture}}}
\newcommand{\evvb}{\raisebox{-0.6cm}{
\begin{picture}(48,37)
\put(0,0){\psbox[scale=0.17]{EVV}}
\end{picture}}}
\newcommand{\wvvb}{\raisebox{-0.6cm}{
\begin{picture}(48,37)
\put(0,0){\psbox[scale=0.17]{WVV}}
\end{picture}}}
\newcommand{\twobox}{\raisebox{-0.5cm}{
\begin{picture}(29,36)
\put(0,0){\psbox[scale=0.2]{fourbox0}}
\end{picture}}}
\newcommand{\twobbox}{\raisebox{-0.5cm}{
\begin{picture}(29,36)
\put(0,0){\psbox[scale=0.2]{fourbox00}}
\end{picture}}}
\newcommand{\twolines}{\raisebox{-0.5cm}{
\begin{picture}(13,36)
\put(0,0){\psbox[scale=0.2]{fourbox2}}
\end{picture}}}
\newcommand{\twocross}{\raisebox{-0.5cm}{
\begin{picture}(13,36)
\put(0,0){\psbox[scale=0.2]{fourbox3}}
\end{picture}}}
\newcommand{\tribox}{\raisebox{-0.5cm}{
\begin{picture}(29,36)
\put(0,0){\psbox[scale=0.2]{sixbox0}}
\end{picture}}}
\newcommand{\tribbox}{\raisebox{-0.5cm}{
\begin{picture}(29,36)
\put(0,0){\psbox[scale=0.2]{sixbox00}}
\end{picture}}}
\newcommand{\threelines}{\raisebox{-0.5cm}{
\begin{picture}(17,36)
\put(0,0){\psbox[scale=0.2]{sixbox1}}
\end{picture}}}
\newcommand{\threecross}{\raisebox{-0.5cm}{
\begin{picture}(19,36)
\put(0,0){\psbox[scale=0.2]{sixbox10}}
\end{picture}}}
\newcommand{\threercross}{\raisebox{-0.5cm}{
\begin{picture}(17,36)
\put(0,0){\psbox[scale=0.2]{sixbox6}}
\end{picture}}}
\newcommand{\threelcross}{\raisebox{-0.5cm}{
\begin{picture}(17,36)
\put(0,0){\psbox[scale=0.2]{sixbox7}}
\end{picture}}}
\newcommand{\threercrosses}{\raisebox{-0.5cm}{
\begin{picture}(17,36)
\put(0,0){\psbox[scale=0.2]{sixbox8}}
\end{picture}}}
\newcommand{\threelcrosses}{\raisebox{-0.5cm}{
\begin{picture}(17,36)
\put(0,0){\psbox[scale=0.2]{sixbox9}}
\end{picture}}}
\newcommand{\nbox}{\raisebox{-0.85cm}{
\begin{picture}(28,50)
\put(0,0){\psbox[scale=0.2]{nbox}}
\end{picture}}}
\newcommand{\nbbox}{\raisebox{-0.85cm}{
\begin{picture}(28,50)
\put(0,0){\psbox[scale=0.2]{nbbox}}
\end{picture}}}
\newcommand{\del}{\raisebox{-0.9cm}{
\begin{picture}(4,52)
\put(0,0){\psbox[scale=0.2]{delta}}
\end{picture}}}
\newcommand{\ijkbox}{\raisebox{-1.1cm}{
\begin{picture}(27,68)
\put(0,0){\psbox[scale=0.25]{tribox}}
\end{picture}}}
\newcommand{\ijkboxx}{\raisebox{-1.3cm}{
\begin{picture}(28,78)
\put(0,0){\psbox[scale=0.25]{triboxx}}
\end{picture}}}
\newcommand{\scalar}{\raisebox{-2.0cm}{
\begin{picture}(214,117)
\put(0,0){\psbox[scale=0.4]{scalar}}
\end{picture}}}
\newcommand{\closedscalar}{\raisebox{-2.5cm}{
\begin{picture}(200,143)
\put(0,0){\psbox[scale=0.4]{closedscalar}}
\end{picture}}}
\newcommand{\spinora}{\raisebox{-0.93cm}{
\begin{picture}(17,53)
\put(0,0){\psbox[scale=0.2]{spinor1}}
\end{picture}}}
\newcommand{\spinorb}{\raisebox{-0.93cm}{
\begin{picture}(17,53)
\put(0,0){\psbox[scale=0.2]{spinor2}}
\end{picture}}}
\newcommand{\spinorc}{\raisebox{-0.93cm}{
\begin{picture}(17,53)
\put(0,0){\psbox[scale=0.2]{spinor3}}
\end{picture}}}
\begin{document}
\title {
\begin{flushright}{\small UT-Komaba 94-8}
\end{flushright}
{\bf\Large Graph Invariants of Vassiliev Type\\ and \\Application to 4D
Quantum
Gravity}}
\author{Nobuharu Hayashi
\thanks{the fellow of the Japan Society for the Promotion of
Science for Japanese Junior Scientists}
\thanks{
{\it E-mail address\/}: hayashin{\char'100}kekvax.kek.jp}
\\
Institute of Physics\\
University of Tokyo\\
Komaba, Tokyo 153, Japan\thanks{This is not the address after April 1.
Any contact to the author is available only via e-mail. }
\\}
\maketitle
%
%
%
%\size{12}{25pt}\selectfont
%
\begin{abstract}
We consider a special class of Kauffman's graph
invariants of rigid vertex isotopy (graph invariants of Vassiliev type).
They are given by a functor from a
category of colored and oriented graphs embedded into a 3-space
to a category of representations of the quasi-triangular ribbon Hopf algebra
$U_q(sl(2,\bf C))$.
Coefficients in  expansions  of them with respect to $x$
($q=e^x$) are known as the Vassiliev invariants of finite type.
In the present paper, we construct two types of tangle operators  of vertices.
One of them corresponds to
a Casimir operator insertion at a transverse double point of Wilson loops.
This paper proposes a non-perturbative generalization of Kauffman's recent
result based on a perturbative analysis of the Chern-Simons
quantum field theory.
As a result, a quantum group analog of
Penrose's spin network is established taking into account of the orientation.
We also deal with the 4-dimensional canonical
quantum gravity of Ashtekar.
It is verified that the graph invariants of Vassiliev type are compatible
with constraints of the quantum gravity in the loop space representation
of Rovelli and Smolin.
\end{abstract}
\pagebreak
\section{Introduction}
The concept of Vassiliev Invariants was originally introduced in the theory
of the knot space from a side of the algebraic topology\cite{VASS1}.Let $M$ be
a space of all smooth  maps $S^1\rightarrow {\bf S}^3$
going
through a fixed point and having a fixed tangent vector at the fixed
point.
$M$ is connected.
The knot space is defined by $M\backslash\Sigma$
where $\Sigma$ is called the "discriminant",
i.e., a set of all singular maps with
multiple points or vanishing tangent vectors.
Any equivalence class of knot embeddings by the ambient isotopy
of $S^3$ corresponds to a connected component of
$M\backslash\Sigma$.
In this sense, all the knots are classified by the space
$H_0(M\backslash\Sigma)$.
Each connected component of the knot space is separated
by walls which constitute the discriminant $\Sigma$.
\par
The Vassiliev invariants are introduced as follows.
Let $\Sigma_j\subset\Sigma$ be a space of $j$-embeddings, i.e., a
space of
all singular maps whose singularities are j transverse double
points.
Then there exists a natural filtration
$\Sigma\supset\Sigma_1\supset\Sigma_2\supset\cdots\supset
\Sigma_j\supset\cdots$.
Using this stratification,
Vassiliev considered a spectral sequence to compute
$H^0(M\backslash\Sigma)$.
He obtained an inclusion:
$H^0(M\backslash\Sigma)\supset \cdots
E^{-j,j}_{\infty}\supset E^{-(j-1),j-1}_{\infty}\cdots
\supset E^{-1,1}_{\infty}$.
Each stabilized limit $E^{-j,j}_\infty$
is called the Vassiliev invariant of order j.
It means that any element of $E^{-j,j}_\infty$ vanishes whenever
it is evaluated for singular maps with more than j transverse double  points.
In addition, Vassiliev conjectured
$E_\infty^{-j,j}\otimes_Z K\cong E_1^{-j,j}\otimes_Z K$
for some rings $K$.
After the pioneered work of Vassiliev,
Birman and Lin\cite{BL1} calculated $E_1^{-j,j}$ in a combinatorial
way based on some axioms.
In their discussions,
they introduced a special class of functionals defined over a space
of
cord diagrams $D^j$.
$D^j$ represents a space of all configurations of j pairs of $2j$
points on $S^1$ connected by j cords.
Inspired by the works of Vassiliev, Birman and Lin,
Kontsevich\cite{Ko1} verified  $E_\infty^{-j,j}\otimes_Z C\cong
E_1^{-j,j}\otimes_Z C$ using an
integral
representation of the Vassiliev invariants.
In the definition of Kontsevich's integral representation,
the weight system\cite{BN1}plays an important role.
It is given by a space of maps from
some kind of Hopf algebra\cite{Ko1} generated by the cord
diagrams
to a ring $K$.
\par
The relation of the Vassiliev invariants to well-known quantum
group
invariants such as the HOMFLY polynomials and the Kauffman
polynomials was investigated
by Birman, Lin\cite{BL1}\cite{Lin1}\cite{Lin2} and
Puinikhin\cite{PU1}\cite{PU2} to
a great extent.
It was pointed out by them that the Vassiliev invariants of finite type or
the Vassiliev invariants of finite order play an important role.
They appear as coefficients in expansions of the quantum group
invariants with respect to $x$where  $q=e^x$.
How about the connection of the Vassiliev invariants to the
Chern-Simons gauge field theory?
In this respect, we can't forget to mention the work of
E. Witten\cite{WIT1}.
He shed a light on a relation between
the quantum group invariants of links
and the 2D conformal field theory
and arrived at an observation that
the quantum group invariants of links can be regarded as vacuum expectation
values of
linking Wilson loops in the CS (Chern-Simons) quantum gauge field
theory.
The Wilson loop is a trace of a holonomy along a link component.
But his discussion was implicitly based on a premise that
two different quantization methods of the CS gauge field theory, i.e.,
the canonical quantization and the path integral quantization,
are equivalent.
Some authors attempted to have a foundation to confirm
themselves about the premise.
A perturbative analysis of the path integral of the CS gauge field
theory is
one of ways to check it.
Works of
Axelrod, Singer\cite{AS1}\cite{AS2} and Bar-Natan\cite{BN1}
treating
the CS perturbation in the context of the Vassiliev invariants are of interest
to us.
In particular,
Bar-Natan uncovered a relation between the Vassiliev Invariants
and the
Feynman diagrams in the CS quantum gauge field theory.
They are connected by the concept of the weight system.
\par
In the present paper,
we are interested in a recent work of Kauffman\cite{KAUFF3}.
(The Vassiliev invariants given by the quantum group invariants of links can be
regarded as
a special class of graph invariants that
he introduced a few years
ago\cite{KAUFF2}\cite{KAUFF3}\cite{KAUFF4}.)
In \cite{KAUFF3}, he discussed
the perturbation theory of the CS quantum gauge field theory to
obtain a path integral representation of the graph invariants of
Vassiliev
type in terms of the Wilson loops with transverse double points.
It was observed by him that the Vassiliev Invariant given by  a j-embedding
can be represented by a CS vacuum expectation  value of the Wilson loops
with j transverse double points where  the quadratic Casimir operators
are inserted.
Such a point of view of the Vassiliev invariants is expected to
provide a neat perspective on the result of Bar-Natan\cite{BN1}.
Another importance of such a point of view comes from an
application to the
4-dimensional non-perturbative canonical quantum gravity of
Ashtekar\cite{AS1}\cite{AS2}.
Kauffman's path integral representation in terms of the Casimir
operator insertion plays
an indispensable role to investigation of the 4D quantum gravity.
\par
This paper is organized as follows.
\S2 contains a brief review on Kauffman's graph
invariants\cite{KAUFF1}\cite{KAUFF2}\cite{KAUFF4} of
rigid vertex isotopy.
We deal with graphs with 4-valent rigid vertices.
The aim of this section is to introduce
Kauffman's graph invariants of Vassiliev type\cite{KAUFF3}.
In \S3, we discuss a generalization of Kauffman brackets of
links to those of singular links whose singularities are
transverse double points. The Kauffman brackets of singular
links are considered provided that
they can be identified with some CS vacuum expectation values of
Wilson loops. In addition to the skein relation,
we introduce the spinor identity, and
explain how it is used to resolve transverse double points.
It is an essential role of the spinor identity
that it allows us to express the Kauffman brackets of singular links
by the Kauffman brackets of links.
In \S4, the discussions
of \S3 are justified in a
context
of representations of the quasi-triangular ribbon Hopf algebra
$U_q(sl(2,\bf C))$.
We attempt to characterize the Kauffman brackets of
singular links in \S3 by
tangle operators of rigid vertex graphs .
They are perceived as a graph generalization of the tangle operators of
links that was considered in the theory of quantum group invariants of
3-manifolds\cite{KM}\cite{RT1}\cite{RT2}.
Our main result is that the graph invariants of Vassiliev type can be expressed
by
two types of tangle operators of rigid vertices which are identified with
physical objects.
We check that our result is compatible with an approximate expression of
Kauffman's graph invariants of Vassiliev type
based on a perturbative analysis of the CS quantum gauge field theory.
In addition to these, we define graph invariants given by graphs composed of
4- and 6-valent vertices
and attempt to find the Casimir insertion representation.
The last section is devoted to physical application to
the 4-dimensional quantum gravity of Ashtekar.
We apply the canonical quantization and consider the loop space representation
of
physical wave functions given by the graph invariants.
They are called spin-network states.
It is verified that all the constraints of the quantum
gravity with vanishing cosmological constant to specify the physical
states are satisfied by the spin network
states given by the graph invariants of Vassiliev type.
This is a subsidiary result of the present paper.
This paper is an enlarged version of \cite{Hayashin} by the physical
application.
\hfill\break
\section{Rigid Vertex Graphs and Graph Invariants of Vassiliev
Type}%%%%%%%%%%%%
\par
This section is devoted to a brief review on Kauffman's graph
invariants that can be extended by the quantum group invariants of links.
Here, we only deals with oriented graphs $G$ whose
constituents are 4-valent rigid vertices and edges and
loops\cite{KAUFF1}\cite{KAUFF2}.
The 4-valent rigid vertex is a disk having four strings emanating
from it.
We shall allow a case in which $G$ is a sum of a finite number of
connected components, i.e., $G=\amalg_i G_i$. We can regard a graph $G$
embedded into
a 3-space $M^3$ as a generalization of a link.
Two embeddings of a graph with 4-valent rigid vertices are said
to
be equivalent, if one is identified with the other by a rigid vertex isotopy of
$M^3$.
In the followings, we introduce the definition of the rigid
vertex isotopy and graph invariants of such a kind of isotopy.
\begin{defn}
An embedding $\phi\colon G\rightarrow M^3$ is called the {\bf
rigid vertex
embedding} if for each vertex $v\in G$, a proper 2-Disk in a ball
neighborhood
B of $\phi(v)$ is specified such that the image of a neighborhood
of $v$ in $G$ is contained in D.
An isotopy $h_t\colon M^3\rightarrow M^3$ between two rigid
vertex
embeddings $\phi_0$ and $\phi_1$ is called the {\bf rigid vertex isotopy} if it
carries through the ball-disk pair for each vertex of $G$.
\end{defn}
For simplicity, we proceed assuming $M^3=S^3$ in the followings.
\par
%%%%%%%%%fig1%%%%%%%%%%%%%%%%%%%%%%%
\begin{figure}
\begin{center}
\epsfile{file=Reidemeister,scale=0.35}
\end{center}
\caption{Generalized Reidemeister Moves generating the rigid
vertex isotopy}
\label{Reidemeister}
\end{figure}
%%%%%%%%%%%%%%%%%%%%%%%%%%%%%%%%%%%
The rigid vertex isotopy is generated by generalized Reidemeister
moves
(see Fig.\ref{Reidemeister}).
The former three moves $R1$, $R2$ and $R3$ are well-known in
the theory of links.
The latter two moves $R4$ and $R5$
\footnote{$R5$ move must be replaced by another one for the
Reidemeister moves generating the topological vertex isotopy\cite{KAUFF1}.
The topological vertex is a point with strings emanating from it.}are
additional ones which  appear in the theory of graphs with
the rigid vertices.
Using a fact that link invariants are defined to be invariant under
the  Reidemeister moves $R1$, $R2$ and $R3$,
we can introduce a huge number of classes of graph invariants
of the rigid vertex isotopy in terms of resolutions of the rigid
vertex.
Kauffman defined such a kind  of graph invariants
in~\cite{KAUFF1}\cite{KAUFF2}.
\par
Suppose that we are given link invariants $P(L)$.
Then, $P(L)$ can always be extended to graph invariants
of the generalized Reidemeister moves.
To be more precise, let's put a resolution of
any one of rigid vertices of the following form
\footnote{Assume that elementary subdiagrams $D$ appearing in $L_D$ have the
blackboard framing.
Furthermore, all the singular links look same outside the
elementary subdiagrams.}:
\begin{equation}
P\left(L^{(j)}\scrooss\right)
=a P\left(L^{(j-1)}\spcross\right)
+b P\left(L^{(j-1)}\sncross\right)
                +c P\left(L^{(j-1)}\sunfold\right)~~~,
\label{res1}
\end{equation}
where $L^{(j)}$ represents an embedding of a graph with j 4-valent
rigid vertices.
Thus $P\left(L^{(j)}\scrooss\right)$
can be expressed by a sum of link invariants after
the resolution of all the vertices.
Then, we find
$
P(G)=\sum_{L\in S}a^{p(S)}b^{n(S)}c^{u(S)}P(L)
$
where $S$ represents a set of $2^j$ links that we obtain after
the resolution of all the rigid vertices of $G$ and
$p(S)$, $n(S)$ and $u(S)$ stand for a number of positive crossings, that of
negative crossings and that of unfoldings respectively.
\begin{thm}
Every $P\left(L^{(j)}\scrooss\right)$ defined by (\ref{res1})
is invariant under the generalized Reidemeister moves.
\end{thm}
{\it Proof.~}
It is immediately proved by definition (see \cite{KAUFF1}\cite{KAUFF2}).
\hfill\qed \vspace{0.2cm}
\break
\par
We are ready to define graph invariants of Vassiliev
type\cite{KAUFF3}.
They are given by putting $a=1$, $b=-1$, and $c=0$, i.e.,
$P\left(L^{(j)}\scrooss\right)
=P\left(L^{(j-1)}\spcross\right)-P\left(L^{(j-1)}\sncross\right)$.
They are objects that we will discuss in the following sections.
%
%
%
%
%%%%%%%%%%%%%%%%%%%%%%%%%%%%%%%%%%%%%%%%%%%%%%%%%%%%%%%%%%%%%
%%%%%%%%%%%%%%%%%%%%%%%%%%%%%%%%%%%%%%%%%%%%%%%%%%%%%%%%%%%%%
\section{Spinor Identity and Kauffman Brackets}%%%%%%%%%%%%
\begin{figure}
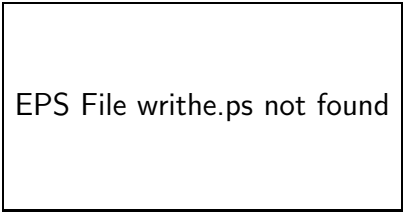

\begin{center}
\epsfile{file=writhe,scale=0.4}
\end{center}
\caption{Signatures for a positive crossing, a negative crossing and a vertex}
\label{writhe}
\end{figure}

%%%%%%%%%%%%%%%%%%%%%%%%%%%%%%%%%%%%%%%%%%%%%%%%%%%%%%%%%%%%%
%%%%%%%%%%%%%%%%%%%%%%%%%%%%%%%%%%%%%%%%%%%%%%%%%%%%%%%%%%%%%%%
The Jones polynomial has a physical analogue in
the CS quantum gauge field theory with $SU(2)$ gauge group.
They are defined to be invariant under the Reidemeister moves
$R1$, $R2$  and $R3$.
Let $P(L)$ be the Jones polynomial of a link $L$ and $Z(L)$
the Kauffman bracket of $L$.
$Z(L)$ is given by a vacuum expectation value of linking
Wilson loops in
the CS  quantum field theory\cite{WIT1}.
The Jones polynomial differs from the Kauffman bracket
only by a phase factor.
To be precise,
$P(L)$ is related to $Z(L)$ by $P(L)=\alpha^{-w(L)}Z(L)$.
$\omega(L)$ is the writhe function given by
a sum of crossing signs of the link diagram of $L$.
It namely means $\omega(L)=\sum_{p\in\Gamma(L)}\epsilon(p)$
where
$\Gamma(L)$ represents a set of all crossings of the link $L$ and a
signature $\epsilon$ takes $1$ ($-1$) for a positive
(negative) crossing (see Fig.\ref{writhe}).
%%%%%%%%%%%%%%%%%%%%%%%%%%%%%%%%%%%%%
%%%%%%%%%%%%%%%%%%%%%%%%%%%%%%%%%%%%%
According to Kauffman's notation
$e(x)=exp(\frac{-i\pi}{k}x)$\cite{KAUFF3},
$\alpha$ is set to be $e(\frac{3}{2})$. Kauffman's brackets are
computable in terms of the skein relation:
\begin{equation}
e(\tfrac1{2})Z\Lpcrossb-e(-\tfrac1{2})Z\Lncrossb=
\left(e(1)-e(-1)\right)Z\Lunfoldingb~~~,
\label{Skein}
\end{equation}
and
\begin{equation}
Z\Lpgraph=e(\tfrac{3}{2})Z\Lgraph,~~~Z\Lngraph=
e(-\tfrac{3}{2})Z\Lgraph~~~.
\label{unfolds}
\end{equation}
Eq.(\ref{unfolds}) gives the reason why
$P(L)$ needs the phase factor $\alpha^{-w(L)}$.
Owing to it, $P(L)$ can be invariant under the Reidemeister move $R1$.
We remark that
the Kauffman brackets are invariants of colored framed links
(colored ribbons).
Indeed,
in the perturbation theory of the CS quantum gauge fields,
we can observe that $Z(L)$ is subject to (\ref{unfolds}) assuming
that
all the Wilson loops are in the 2-dimensional representation of
$SU(2)$. The Jones polynomials are also computable in terms of
the skein relation obtained by rewriting (\ref{Skein}) :
\begin{equation}
e(2)P\Lpcrossb-e(-2)P\Lncrossb=\left(e(1)-e(-1)\right)P\Lunfoldingb~~~.
\label{Skein2}
\end{equation}
We conventionally put an initial condition
$P\Lcircleb=1$ for an unknot.
\par
In the followings, assume that we can identify the Kauffman brackets
$Z(L)$ with the CS vacuum expectation values of Wilson
loops.
The Wilson loop is given by
a trace of a holonomy along any one of link components $L_i$
where $L=\amalg_i L_i$.
We denote it by
$W^{(i)}(L_i,A)\equiv Tr\left(U(s,s)\right)=Tr\left(Pexp\left(i\oint_{\gamma_i}
A\right)\right)$.
$P$ represents the path-ordered product along a closed path
$\gamma_i(s)$ ($0\leq s\leq 1$) which describes the embedding of
the link component $L_i$.
$Pexp\left(i\oint_{\gamma_i} A\right)$ is the holonomy along $L_i$
given by a connection
1-form of a $SU(2)$-principle bundle over $M^3$.
\par
Let's discuss a generalization of the
Kauffman brackets of links to those of singular links whose
singularities are transverse double points.
First, we introduce a vacuum
expectation value of  N linking Wilson loops
with only one transverse double point:
\begin{equation}
\begin{split}
Z\left(L^{(1)}\scross\right)
     &\equiv\int DA
       Tr\left(U_j(s,t)U_j(t,s)\right)\cdot
\Pi_{i\ne j}^{N}W^{(i)}(L_i,A)\\
                    &\phantom{=}
            \times exp\left(\frac{ik}{4\pi}\int
                     Tr(A\wedge dA
+\frac{2}{3}A\wedge A\wedge A)\right)
\end{split}
\end{equation}
for case 1 in which the j-th link component transversely intersects only at
one point ($\gamma_j(s)=\gamma_j(t)$ and $0\leq s<t\leq 1$), and
\begin{equation}
\begin{split}
Z\left(L^{(1)}\scross\right)
     &\equiv\int DA
       Tr\left(U_j(s,s)\right)\cdot Tr\left(U_k(t,t)\right)\cdot
           \Pi_{i\ne j,k}^{N}W^{(i)}(L_i,A)\\
                    &\phantom{=}
            \times exp\left(\frac{ik}{4\pi}\int
                     Tr(A\wedge dA
                    +\frac{2}{3}A\wedge A\wedge A)\right)
\end{split}
\end{equation}
for case 2 in which the j-th and k-th link components transversely
intersect only at one point.
$\gamma_j(s)$ and $\gamma_k(t)$ describe $L_j$ and $L_k$
respectively and
satisfy $\gamma_j(s)=\gamma_k(t)$.
Another important object that we need to
introduce for the latter  argument is $Z\Lonecasimirb$.
Let's define it as follows.
%%%%%%%%%%%%%%%%%%%%%%%%%%%%%%%%%%%%%%%%%%%%%%%
\begin{figure}
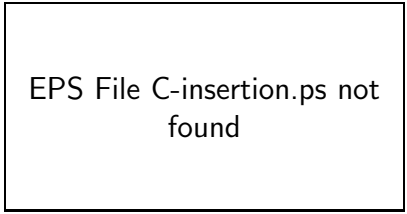

\begin{center}
\epsfile{file=C-insertion,scale=0.5}
\end{center}
\caption{Casimir Insertion at a transverse double point}
\label{Casimir}
\end{figure}
%%%%%%%%%%%%%%%%%%%%%%%%%%%%%%%%%%%%%%%%%%%%%%%
\begin{defn}
Operation of inserting the quadratic Casimir
operator at a transverse double point
is called the {\bf Casimir insertion} (see Fig.\ref{Casimir}).
Suppose that N Wilson loops have only one transverse double point.
The CS vacuum expectation value of N Wilson loops with the Casimir
insertion is given  by
\begin{equation}
\begin{split}
Z\left(L^{(1)}\sccross\right)
     &\equiv\int DA\sum_{a=1}^{dim(SU(2))=3}
       Tr\left(U_j(s,t)T_aU_j(t,s)T_a\right)\cdot
           \Pi_{i\ne j}^{N}W^{(i)}(L_i,A)\\
                    &\phantom{=}
            \times exp\left(\frac{ik}{4\pi}\int
                     Tr(A\wedge dA
                    +\frac{2}{3}A\wedge A\wedge A)\right)
\end{split}
\end{equation}
for case 1,and
\begin{equation}
\begin{split}
Z\left(L^{(1)}\sccross\right)
     &\equiv\int DA\sum_{a=1}^{dim(SU(2))=3}
       Tr\left(U_j(s,s)T_a\right)\cdot Tr\left(T_aU_k(t,t)\right)\cdot
           \Pi_{i\ne j,k}^{N}W^{(i)}(L_i,A)\\
                    &\phantom{=}
            \times exp\left(\frac{ik}{4\pi}\int
                     Tr(A\wedge dA
                    +\frac{2}{3}A\wedge A\wedge A)\right)
\end{split}
\end{equation}
for case 2.
\end{defn}
A generalization to more complicated cases in which there exist more than one
transverse
double point is trivial.
In such more general cases,
it is convenient to introduce the following notation.
Let $Z(L^{(k,j-k)})$ ($0\leq k\leq j$) be generalized Kauffman
brackets given by a singular link whose singularities are j transverse double
points.
The Casimir operators are inserted at $j-k$ transverse double points one by
one.
In the notation,
$Z(L^{(0,0)})$ is nothing but the Kauffman bracket $Z(L)$ of a link $L$.
$Z\left(L^{(1,0)}\right)$ and
$Z\left(L^{(0,1)}\right)$ correspond to
$Z\left(L^{(1)}\scross\right)$ and
$Z\left(L^{(1)}\sccross\right)$ respectively.
Some of the generalized Kauffman brackets
$Z\left(L^{(k,j-k)}\right)$ are
not independent  because of the Fierz identity, i.e.,
$\sum_a(T_a)_{ij}(T_a)_{kl}=\frac{1}{2}\delta_{il}\delta_{jk}
                           -\frac{1}{2N}\delta_{ij}\delta_{kl}$
($N=2$ for $SU(2)$).
It induces a relation among them:
\begin{equation}
Z\left(L^{(k,j-k)}\sccross\right)
=\frac{1}{2}Z\left(L^{(k,j-k-1)}\sunfold\right)
                     -\frac{1}{4}
                     Z\left(L^{(k+1,j-k-1)}\scross\right)~~~.
\label{Fierz}
\end{equation}
\par
Regarding the generalized Kauffman brackets as
vacuum expectation values of intersecting Wilson loops,
we must take into account of the spinor identity
in addition to the skein relation (\ref{Skein}).
It is satisfied by $Z\left(L^{(j,0)}\right)$ for $j\geq 1$.
If there exist some Casimir insertions,
we can always eliminate them using the Fierz identity (\ref{Fierz}).
Let's explain the spinor identity below.
Let $A$ and $B$ be elements of $SU(2)$, i.e.,
invertible $2\times 2$ matrices.  It follows that
$Tr(A)Tr(B)=Tr(AB)             +Tr(AB^{-1})$. This is derived from a fact :
$\epsilon_{ab}\epsilon^{cd}
                =\delta_a^c\delta_b^d-\delta_a^d\delta_b^c$.
The spinor identity induces a relation among the
generalized Kauffman brackets:
\begin{equation}
Z(\alpha\cup\beta)=Z(\alpha\beta)+Z(\alpha\beta^{-1})~~~.
\label{spinorKauff}
\end{equation}
On the left hand side, $Z(\alpha\cup\beta)$ represents a vacuum
expectation
value of a product of two Wilson loops
$Tr\left(U(\alpha)\right)Tr\left(U(\beta)\right)$ given by two
closed paths
$\alpha$ and $\beta$ with the common base point.
On the right hand side,
$Z(\alpha\beta)$ ($Z(\alpha\beta^{-1})$) represents a vacuum
expectation
value of a Wilson loop $Tr\left(U(\alpha\beta)\right)$
($Tr\left(U(\alpha\beta^{-1})\right)$) along a composite and closed
path
$\alpha\beta$  ($\alpha\beta^{-1}$).
(\ref{spinorKauff}) is called the spinor identity of the generalized Kauffman
brackets.
\par
\begin{figure}
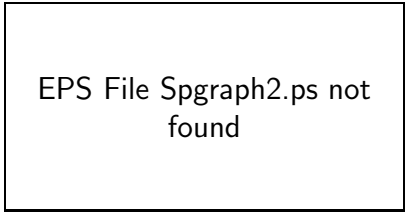

\begin{center}
\epsfile{file=Spgraph2,scale=0.4}
\end{center}
\caption{In the case 1,
the composite path $\alpha\beta$ makes a transverse double point
at the base point. In the case 2, closed paths $\alpha$ and $\beta$
transversely intersect at the base point.
}
\label{Spgraph2}
\end{figure}

We can use the spinor identity to resolve transverse double
points.
It is enough to consider two cases (see Fig.\ref{Spgraph2}).
%%%%%%figure%%%%%%%%%%%%%%%%%%%%%%
%%%%%%%%%%%%%%%%%%%%%%%%%%%%%%%%%%
In one case (case 1), the closed paths $\alpha$ and $\beta$ compose a
closed path $\alpha\beta$ so that
a transverse double point appears at their base point.
Then the spinor identity takes the following form:
\footnote{We employ a notation
$Z(D)$ instead of $Z(L^{(j)}_D)$ below.
$D$ represents an
elementary sub-diagram of
a singular link $L^{(j)}$.
Any diagram depicted in the followings possesses
a generalized blackboard framing.
In this framing,
any transverse double point should be replaced by a disk parallel to
the paper.}
\begin{equation}
Z\Unfoldinb=Z\Crosssb+Z\Unfoldingssb~~~,
\label{spinor1}
\end{equation}
where
\begin{equation*}
Z(\alpha\cup\beta)=Z\Unfoldinb,~~~
Z(\alpha\beta)=Z\Crosssb,~~~
Z(\alpha\beta^{-1})=Z\Unfoldingssb~~~.
\end{equation*}
In the other case (case 2), two closed paths $\alpha$ and $\beta$
transversely intersect at their base point.
In this case, the spinor identity is of the form:
\begin{equation}
Z\Crossssb=Z\Unfoldinnb+Z\Unfoldingsb~~~,
\label{spinor2}
\end{equation}
where
\begin{equation*}
Z(\alpha\cup\beta)=Z\Crossssb,~~~
Z(\alpha\beta)=Z\Unfoldinnb,~~~
Z(\alpha\beta^{-1})=Z\Unfoldingsb~~~.
\end{equation*}
We should remark that a transverse double point appears in
$\alpha\beta$ in the former case, and in $\alpha\cup\beta$ in
the latter case.
\par
\begin{figure}
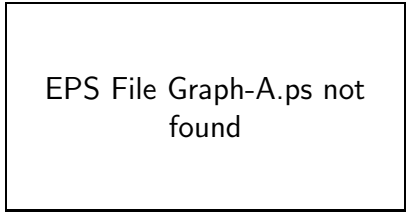

\begin{center}
\epsfile{file=Graph-A,scale=0.3}
\end{center}
\caption{In the diagram $G_a$, $\alpha$ and $\beta$ are closed
paths.
The composite path $\alpha\beta$ makes a transverse double point.
In the diagram $G_b$, circles
$\alpha$ and $\beta$ transversely intersect only at their base
point.
}
\label{Graph-A}
\end{figure}

Let's prepare a few results derived from the
spinor
identity necessary in the latter argument.
First, let's discuss such a diagram as $G_a$ (see Fig.\ref{Graph-A}).
Suppose that a closed path $\beta$ has no intersection point except
at the base point.
%%%%%%figure Graph-A%%%%%%%%%%%%%%%%%%%%%%
%%%%%%%%%%%%%%%%%%%%%%%%%%%%%%%%%%
Then it follows that
\begin{equation}
Z(\alpha\cup\beta)=Z\Graphab=
(e(1)+e(-1))Z\Graphac
\label{circl}
\end{equation}
from the skein relation (\ref{Skein}).
Noticing $Z(\alpha\beta^{-1})=Z\Graphac$ and applying the spinor
identity (\ref{spinor1}) to computation of $Z(\alpha\beta)$, we
immediately obtain
\begin{equation}
Z(\alpha\beta)=Z\Graphaa=(e(1)+e(-1)-1)Z\Graphac~~~.
\label{cups1}
\end{equation}
Second, we consider such a diagram as $G_b$ (see
Fig.\ref{Graph-A}) where two circles $\alpha$ and $\beta$
transversely
intersect at their base point.
Choosing the normalization $Z\Circleb=1$,
it is obvious that
$Z(\alpha\cup\beta)=Z\Graphba$,
$Z(\alpha\beta)=Z\Graphbd=e(\frac{3}{2})$ and
$Z(\alpha\beta^{-1})=Z\Graphbe=e(-\frac{3}{2})$.
We immediately find substituting them into (\ref{spinor2})
\begin{equation}
Z\Graphba=Z\Graphbd+Z\Graphbe=
         e(\tfrac{3}{2})+e(-\tfrac{3}{2})~~~.
\end{equation}
The two cases play an important role in proving
the following proposition.
\begin{prop}
Let's assume that
$Z\left(L^{(1,0)}\scross\right)$ is uniquely expressible as a  sum of
$Z\left(L^{(0,0)}\spcross\right)$ and
$Z\left(L^{(0,0)}\sncross\right)$.
Then it must hold that
\begin{equation}
Z\left(L^{(1,0)}\scross\right)
        =\frac{1}{e(\frac{1}{2})+e(-\frac{1}{2})}
\left(Z\left(L^{(0,0)}\spcross\right)+
Z\left(L^{(0,0)}\sncross\right)\right)~~~.
\label{resporp}
\end{equation}
\end{prop}
\par
{\it Proof~~~}
By assumption, we can put
$Z\left(L^{(1,0)}\scross\right)
=c_1Z\left(L^{(0,0)}\spcross\right)
+c_2Z\left(L^{(0,0)}\sncross\right)$.
For simplicity, taking into account of the Fierz identity (\ref{Fierz}),
it is convenient to put instead of it
{\allowdisplaybreaks
\begin{align}
e(\tfrac{1}{2})Z\left(L^{(0,0)}\spcross\right)
+e(-\tfrac{1}{2})Z\left(L^{(0,0)}\sncross\right)
&=a_1Z\left(L^{(0,1)}\sccross\right)
+a_2Z\left(L^{(1,0)}\scross\right)\notag\\
&=\left(a_2-
\frac{a_1}{4}\right)Z\left(L^{(1,0)}\scross\right)
+\frac{a_1}{2}Z\left(L^{(0,0)}\sunfold\right)~~~.
\label{Pa}
\end{align}}
The determination of the coefficients $a_1$ and $a_2$ finishes the proof the
proposition 3.1.
\par
First, to determine the coefficients $a_1$ and $a_2$,
let's consider a diagram like $G_a$ depicted in Fig.\ref{Graph-A}.
Then (\ref{Pa}) becomes
{\allowdisplaybreaks
\begin{align}
&e(\tfrac{1}{2})Z\Graphap+e(-\tfrac{1}{2})Z\Graphan\notag\\
&\phantom{e(\tfrac{1}{2})Z\Graphap}
=\left(a_2-\frac{a_1}{4}\right)Z\Graphaa
            +\frac{a_1}{2}Z\Graphab~~~.
\label{Ga1}
\end{align}}
Remembering (\ref{unfolds}) , (\ref{circl}) and (\ref{cups1}), one can find
a relation between $a_1$ and $a_2$:
{\allowdisplaybreaks
\begin{align}
e(2)+e(-2)&=\frac{a_1}{4}\left(e(1)+e(-1)+1\right)
+a_2\left(e(1)+e(-1)-1\right)~~~,\notag\\
\intertext{or}
a_2&=\frac{e(2)+e(-2)}{e(1)+e(-1)-1}-\frac1{4}
\frac{e(1)+e(-1)+1}{e(1)+e(-1)-1}a_1~~~.
\label{coef}
\end{align}}
\par
What we have to do next is to determine the coefficient $a_1$.
Let's consider the diagram $G_b$ depicted
in Fig.\ref{Graph-A}.
Application of the resolution (\ref{Pa}) to the diagram $G_b$ gives
{\allowdisplaybreaks
\begin{align}
&e(\tfrac1{2})Z\Graphbb+e(-\tfrac1{2})Z\Graphbc\notag\\
&=\left(a_2-
\frac{a_1}{4}\right)Z\Graphba+\frac{a_1}{2}Z\Graphbd~~~.
\label{Hopfa}
\end{align}}
Three of the four terms appearing in (\ref{Hopfa}) are already
computed.
It is also easy to compute the rest term, i.e.,
the first term on the left hand side of (\ref{Hopfa}).
After some algebra, we find
{\allowdisplaybreaks
\begin{align}
Z\Graphbb&=e(2)+e(-2)~,\quad Z\Graphbc=e(1)+e(-1) ~~~,\notag\\
Z\Graphba&=e(\tfrac{3}{2})+e(-\tfrac{3}{2})~,\quad
Z\Graphbd=e(\tfrac{3}{2})~~~.
\label{Hopfb}
\end{align}}
After substitution of (\ref{Hopfb}) into (\ref{Hopfa}),
we finally obtain the solution
\begin{equation}
a_1=2\left(e(1)+e(-1)-2\right)~~~.
\label{prop1}
\end{equation}
\par
Thus the coefficients $a_1$ and $a_2$ are determined.
$Z\left(L^{(1,0)}\scross\right)$ can be expressed in terms
of  $Z\left(L^{(0,0)}\spcross\right)$
and $Z\left(L^{(0,0)}\sncross\right)$ according to (\ref{Pa}) and
the skein relation (\ref{Skein}).
The proposition is proved.
\hfill\qed \vspace{0.2cm}
\break
For the Casimir insertion representation
$Z\left(L^{(0,1)}\sccross\right)$,
it follows from the Fierz
identity (\ref{Fierz}) that
\begin{equation}
Z\left(L^{(0,1)}\sccross\right)=\frac{1}
{4\left(e(\frac{1}{2})-e(-\frac{1}{2})\right)}
\left(Z\left(L^{(0,0)}\spcross\right)
-Z\left(L^{(0,0)}\sncross\right)\right)~~.
\end{equation}
There is no reason of restriction to the $j=1$ case.
It comes from simplicity of the computation.
Indeed, we can consider
more general
cases in which there are more than one transverse
double point.
We can expect the following.
\begin{conj}
Let $Z\left(L^{(j,0)}\right)$ ($j\geq 1$) be the generalized Kauffman
bracket
given by a singular link with j transverse double points.
It is given by resolution of any one of the j transverse double points:
\begin{equation}
Z\LLjcrossb=\frac{1}{e(\frac{1}{2})+e(-\frac{1}{2})}
\left(Z\LLjpcrossb+Z\LLjncrossb\right)~~~.
\label{resporp2}
\end{equation}
\end{conj}
We checked that there is no contradiction for several cases.
The next section justifies the conjecture in a context of
the representations of the quasi-triangular ribbon Hopf algebra
$U_q(sl(2,\bf C))$.
\section{Penrose's Spin-network and Graph Invariants of Vassiliev
Type}%%%%%%%%%%%%%%%%%
\par
\begin{figure}
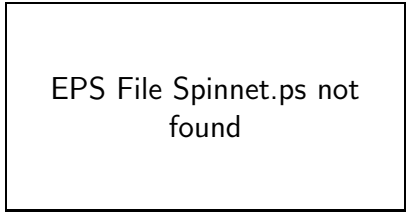

\begin{center}
\epsfile{file=Spinnet,scale=0.45}
\end{center}
\caption{A trivalent graph}
\label{Spinnet}
\end{figure}

Let's begin with a brief introduction
to Penrose's spin-network\cite{Penrose1}\cite{Penrose2}.
It is a trivalent graph whose edges are
in representations of $SU(2)$ (see Fig.\ref{Spinnet}).
Penrose's spin-network has a strand representation in which
each strand is colored by a tensor product representation of
fundamental representations.
%%%%%%figure%%%%%%%%%%%%%%%%%%%%%
%%%%%%%%%%%%%%%%%%%%%%%%%%%%%%%%%%
In the strand representation, each trivalent vertex
is coupling of three strands (Fig.\ref{trivertex}), e.g.,
coupling of i-, j-, k-units is expressed as
\begin{equation}
\ijkbox=\ijkboxx~~~.
\label{trivalent}
\end{equation}
The induces $(i,j,k)$  are admissible only if $i+j+k\geq 2max(i,j,k)$ and
$i+j+k\in
2\bf Z$. Any box stands for an operator composed of skew-symmetrizers and
symmetrizers.
It may be projection operators onto irreducible representations of
$SU(2)$
(the Yang tableau operators\cite{Penrose1}).
To be concrete,
the skew-symmetrizer for n-units is
\begin{equation}
\nbox\equiv\frac{1}{n!}\sum_{\sigma\in {\frak S}_n}sign(\sigma)
\delta^{i_1}_{j_{\sigma(1)}}\delta^{i_2}_{j_{\sigma(2)}}
\cdots\delta^{i_n}_{j_{\sigma(n)}}~~~,~\text{where}~~~\del=\delta^i_j~~~.
\label{nskew}
\end{equation}
The indices $i_m$, $j_m$ ($1\leq m \leq n$) represent spinor indices.
${\frak S}_n$ represents the symmetric group of order n.
For instance, the skew-symmetrizers of two- and three-units are
{\allowdisplaybreaks
\begin{equation}
\frac{1}{2!}\twobox=\twolines-\twocross,~~~~~
\frac{1}{3!}\tribox=
\threelines+\threelcrosses+\threercrosses-\threelcross-\threercross-
\threecross~~~.
\label{skew}
\end{equation}}
The following is the symmetrizer which is given by n-units:
\begin{equation}
\nbbox\equiv\frac{1}{n!}\sum_{\sigma\in {\frak S}_n}
\delta^{i_1}_{j_{\sigma(1)}}\delta^{i_2}_{j_{\sigma(2)}}
\cdots\delta^{i_n}_{j_{\sigma(n)}}~~~~,
\label{nsym}
\end{equation}
The symmetrizers of two- and three-units,
\begin{equation}
\frac{1}{2!}\twobbox=\twolines+\twocross,~~~~~
\frac{1}{3!}\tribbox=\threelines+\threelcrosses+\threercrosses+\threelcross+
\threercross+\threecross~~~.
\label{sym}
\end{equation}
\par
Penrose's spin-network is, furthermore,  specified by
the spinor identity ($SU(2)$ case of Mandelstam's identity
\cite{Giles}),
\begin{equation}
\sqrt{-1}\epsilon^{ab}\sqrt{-1}\epsilon_{cd}
-\delta^a_d\delta^b_c+\delta^a_c\delta^b_d=0~~~\Rightarrow~~~
\spinora-\spinorb+\spinorc=0~~~~.
\label{spinoridenty}
\end{equation}
%%%%%%figure%%%%%%%%%%%%%%%%%%%%%%
\begin{figure}
\begin{center}
\epsfile{file=trivertex,scale=0.4}
\end{center}
\caption{Strand representation of a trivalent vertex}
\label{trivertex}
\end{figure}
%%%%%%%%%%%%%%%%%%%%%%%%%%%%%%%%
\par
Invariants of Penrose's spin-network are given by a functor from the category
of
trivalent graphs to the category of representations of $SU(2)$.
(In general, we can define them for $SL(2,\bf  C)$
where we must take account of primed and unprimed spinors.)
The invariants of Penrose's spin-network
are given by contraction of tensor representations of the spinor
representations.
It corresponds to making scalars in the abstract tensor system
(ATS\cite{Penrose2}).
Let's introduce a type (n,m)-tensor
$\varGamma_{l_1\cdots l_n}^{k_1\cdots k_m}$.
It is a morphism from
$V^{\otimes n}$ (the tensor product of n spinors) to $V^{\otimes m}$, i.e.,
is in ${V^{\otimes n}}^*\otimes V^{\otimes m}$.
Especially, a  morphism from $V^0$ (the trivial representation) to $V^{\otimes
n}$
is here denoted by $\varPsi^{l_1\cdots l_n}$.
A morphism from $V^{\otimes m}$ to $V^0$ is denoted by $\varPhi_{k_1\cdots
k_m}$.
The indices are raised or lowered by the Levi-Civita tensor.
A graphical expression of invariants takes the following form:
{\allowdisplaybreaks
\begin{align}
&\sum_{i,j}
\varPhi_{i_1i_2\cdots i_{2n-1}i_{2n}}
T_{j_1j_2\cdots j_{2n-1}j_{2n}}^{i_1i_2\cdots i_{2n-1}i_{2n}}
\varPsi^{j_1j_2\cdots j_{2n-1}j_{2n}}\notag\\
&=\scalar
{}~~~.
\label{scalar}
\end{align}}
Let's consider a closed spin-network as depicted in the Fig.\ref{Spinnet}.
The invariants are given by contraction of
$T_{j_1j_2\cdots j_{2n-1}j_{2n}}^{i_1i_2\cdots i_{2n-1}i_{2n}}$
in terms of the Levi-Civita tensor.
The graphical expression is
{\allowdisplaybreaks
\begin{align}
&\sum_{i,j}\sqrt{-1}\epsilon_{i_1i_2}\cdots\sqrt{-1}\epsilon_{i_{2n-1}i_{2n}}
T_{j_1j_2\cdots j_{2n-1}j_{2n}}^{i_1i_2\cdots i_{2n-1}i_{2n}}
\sqrt{-1}\epsilon^{j_1j_2}\cdots\sqrt{-1}\epsilon^{j_{2n-1}j_{2n}}\notag\\
&=
\closedscalar~~~.
\label{closedscalar}
\end{align}}
In the theory of Penrose's spin-network,
invariants are regarded as statistical partition functions associated with a
discretized spacetime.
In the continuum limit, they are expected to coincide with
exponential functions of actions in the general relativity.
Such a framework has been referred to as the Regge calculus\cite{Regge}.
\par
The next subsection is devoted to the q-analog of Penrose's spin-network.
We follow the Reshetikhin-Turaev construction of the 3-manifold invariants
based on the quantum groups.
The unit of strands is
in the 2-dimensional representation of $U_q(sl(2,\bf C))$ with q a root of
unity.
In the q-analog of Penrose's spin-network, we must take account of
the spinor identity (\ref{spinoridenty})
in the context of the representation theory of $U_q(sl(2,\bf C))$
(exactly, the quasi-triangular ribbon Hopf algebra given by $U_q(sl(2,\bf
C))$).
We must remember the discussion in \S3
in which the spinor identity can be translated into
the language of the CS quantum gauge field theory.
It is satisfied by vacuum expectation values of intersecting Wilson loops
(see \ref{spinorKauff} in \S 2).
The q-analog of Penrose's spin-network can be interpreted in many ways indeed.
We here mean the q-analog by existence of unit of strands that
satisfies the spinor identity in the CS quantum field theory.
The spinor identity plays a crucial role in defining operators of rigid
vertices of graphs.
The spinor identity and the Fierz identity lead us
to the Casimir insertion representation of the graph invariants of Vassiliev
type.
The q-analog of trivalent coupling of strands is a tedious exercise
in the theory of the quantum group invariants.
Such a discussion is out of the purpose of the present paper.
%
%
%
%
%
%
%
%
%
%
%
%
%%%%%%%%%%%%%%%%%%%%%%%%%%%%%%%%%%%%%%%%%%%%%%%%%%%%%%%%%%%%%%%%%%%%%%%%%%%%
%%%%%%%%%%%%%%%%%%%%%%%%%%%%%%%%%%%%%%%%%%%%%%%%%%%%%%%%%%%%%%%%%%%%%%%%%%
\subsection{Reshetikhin-Traev Construction of Graph Invariants}%%%%%%%%%%%%%
%%%%%%%%%%%%%%%%%%%%%%%%%%%%%%%%%%%%%%%%%%%%%%%%%%%%%%%%%%%%%%%%%%%%%%%%%%%
%%%%%%%%%%%%%%%%%%%%%%%%%%%%%%%%%%%%%%%%%%%%%%%%%%%%%%%%%%%%%%%%%%%%%%%%%%%
According to the work of Reshetikhin and
Turaev\cite{RT1}\cite{RT2}\cite{Turaev},
the tangle operators are given by a functor of the category of
tangles
to the category of representations
of the quasi-triangular ribbon Hopf algebra .
The quasi-triangular ribbon Hopf algebra is specified by
 the Hopf algebra $A$, the universal $R$ matrix and an invertible charmed
element.
The tangle represents an object composed not only of a set of braiding open
strings
but also of pair creation and pair annihilation of them.
The tangle operators of links such as the HOMFLY polynomials
are composed of tangle operators of
elementary tangles of link diagrams, i.e.,
a positive crossing, a negative crossing, a pair creation and a
pair annihilation.
The quantum group invariants of Reshetikhin and Turaev
provided a powerful tool to classify 3-manifolds, i.e., quantum group
invariants of
3-manifolds in terms of the Jones polynomials.
Then Kirby and Mervin made their discussion more precise.
The invariants of 3-manifolds are defined to be invariant under the Kirby-move.
Such invariance comes from a fact that there are no less than one way to obtain
any closed and oriented 3-manifold by the Dehn surgery along a link from $S^3$
{}.
\par
In this sub-section,
we discuss a generalization of tangle operators of
links to those of rigid vertex graphs in the context of
Reshetikhin-Traev construction of topological invariants.
The tangle operators of graphs are also given by a functor from
the category of graphs to the category of representations
of the quasi-triangular ribbon Hopf algebra.
The graph invariants of Vassiliev type introduced in \S2 are
regarded as a special class of the tangle operators of rigid vertex graphs.
Our main interest is whether the graph invariant of Vassiliev type
(given by the tangle operator) can be identified
with any vacuum expectation value of the
Wilson loops with a finite number of double points or not.
In the followings, our discussion is restricted to
the quasi-triangular ribbon Hopf algebra given by
$U_q(sl(2,\bf C))$.
We mainly follow  the work of Kirby and Mervin\cite{KM} in regard of notation
and
definition of the tangles and the tangle operators.
\par
It is well known that
the Jones polynomial given by a link $L$ defers by a phase
factor from the Kauffman bracket.
Assuming that all link
components of $L$  are in the 2-dimensional representation of
$U_q(sl(2,\bf C))$, it follows that
$P(L)=\alpha^{-\omega(L)}F(L)$.
The Kauffman bracket $F(L)$ is given by a functor from the
category of
framed links to the category of representations of the quasi-
triangular Hopf algebra
\footnote{
In the Reshetikhin-Traev construction, $F$ is referred
to as a functor from the category of graphs to the category of
representations of $U_q(sl(2,\bf C))$.
The link is a graph with no legs.
The quantum group invariants of links are given
by the morphism $F:L\rightarrow \bf C$.
} .
On the other hand, remember that the Kauffman
brackets are regarded as the CS vacuum expectation values
of the Wilson loops in \S3 where they are denoted by $Z(L)$.
We must emphasize on a fact that the Kauffman bracket
$F(L)$ is based on the quantum group
and the Kauffman bracket $Z(L)$ is based on the CS path integral.
In the present sub-section, they are distinguished.
\par
The tangle operators of links can be
extended to those of rigid vertex graphs
in terms of resolution of rigid vertices following
Kauffman's graph extension theorem (explained in \S 2).
The graph invariants of Vassiliev
type are given by $P\Lcroossb=P\Lpcrossb-P\Lncrossb$.
Let's introduce operators (Kauffman brackets) $F(G)$ of rigid vertex graphs
such that
$P(G)=\alpha^{-\omega(G)}F(G)$.
\footnote{$\omega(G)=\sum_{p\in\Gamma(G)}\epsilon(p)$. $\Gamma(G)$ represents a
set of
all crossings and rigid vertices of $G$.
The signature for the rigid vertex takes zero.}
We wish to investigate how the operators $F(G)$ are related
to vacuum expectation values of transversely intersecting Wilson loops.
We proceed as follows.
First, we aim at defining
tangle operators of rigid vertices such that
the graph invariants composed of them can be naturally identified with
$Z\left(L^{(k,j-k)}\right)$ (introduced in \S3), i.e.,
the Kauffman brackets given
by the vacuum expectation values of the Wilson loops.
Second, we attempt to express the operators $F(G)$ by them.
\par
Let's begin with defining two types of tangle operators of rigid vertices.
It will be soon shown that operators composed of them have significant
properties that
allow us to identify them with such CS vacuum expectation values of
the Wilson
loops with transverse double points as was discussed in \S3.
\begin{defn}
Let \{$L^{(k,j-k)}$\} ($j\geq1$ and $0\leq k\leq j$)
be a set of graphs with j 4-valent rigid vertices.
Each $L^{(k,j-k)}$ has j-k vertices marked by \copyright.
The rest of k vertices are with no mark.
We introduce tangle operators $F\left(L^{(k,j-k)}\right)$.
They are defined by two types of resolution
of the rigid vertices.
Resolution of a vertex with no mark is given by
\begin{equation}
F\left(L^{(k,j-k)}\scross\right)
=\frac1{e(\tfrac1{2})+e(-\tfrac1{2})}
               \left(F\left(L^{(k-1,j-k)}\spcross\right)+
 F\left(L^{(k-1,j-k)}\sncross\right)\right)~~~,
\label{resolution}
\end{equation}
and resolution of a vertex marked by \copyright is given by
{\allowdisplaybreaks
\begin{align}
F\left(L^{(k,j-k)}\sccross\right)
&=\frac1{2}F\left(L^{(k,j-k-1)}\sunfold\right)
-\frac1{4}F\left(L^{(k+1,j-k-1)}\scross\right)
\label{tanglecasimir}\\
&=\frac1{4\left(e(\tfrac1{2})-e(-\tfrac1{2})\right)}
\left( F\left(L^{(k,j-k-1)}\spcross\right)
               -F\left(L^{(k,j-k-1)}\sncross\right)\right)~~~.\notag
\end{align}}
\end{defn}
It is easily noticed that every tangle operator
$F\left(L^{(k,j-k)}\right)$
can be completely determined by the tangle operators of links
by iteration. Notice that (\ref{tanglecasimir}) is nothing but the Fierz
identity.
The rest of this sub-section will be devoted to proving
that  $F\left(L^{(j,0)}\right)$ satisfies the spinor identity
(see (\ref{spinor1}) and (\ref{spinor2})).
\par
Let's start with the following theorem for $j=1$ case.
\begin{thm}
Let $F\left(L^{(1,0)}\right)$ be the tangle operator of a graph with only one
rigid vertex. It is given by the resolution (\ref{resolution}).
Then it satisfies the spinor identity
\begin{equation}
F\left(L^{(1,0)}\scross\right)=F\left(L^{(0,0)}\sunfold\right)
                  -F\left(L^{(0,0)}\suunfold\right)~~~~,
\label{tanglecross2}
\end{equation}
for case 1 (see the left diagram of Fig.\ref{Spgraph3}) ,
%%%%%%figure%%%%%%%%%%%%%%%%%%%%%%
\begin{figure}
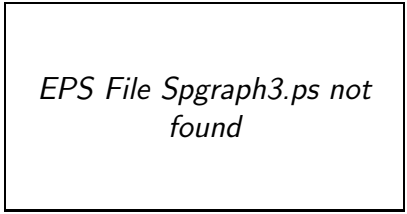

\begin{center}
\epsfile{file=Spgraph3,scale=0.4}
\end{center}
\caption{Two cases of resolution of a rigid vertex concerning the spinor
identity are depicted.
The graph of case 1 is obtained
from the graph of case 2 by interchanging the two
points on the rigid vertex (the disk) from which
the two out-going strings emanate.
These graphs have one-to-one correspondence
to singular links with transverse double points. }
\label{Spgraph3}
\end{figure}
%%%%%%%%%%%%%%%%%%%%%%%%%%%%%%%%%%
and
\begin{equation}
F\left(L^{(1,0)}\scross\right)=
F\left(L^{(0,0)}\sunfold\right)+F\left(L^{(0,0)}\stunfold\right)~~~,
\label{tanglecross1}
\end{equation}
for case 2 (see the right diagram of Fig.\ref{Spgraph3}).
\end{thm}
{}From the second terms on the right hand sides
of (\ref{tanglecross2}) and (\ref{tanglecross1}) ,
the spinor identity is accompanied by the orientation reverse.
Before proving the theorem 4.1, we need the following lemma.
\begin{lem}
Let all the link components of a link $L$ be
in the 2-dimensional representation $V^2$ and
$k\in Z$ be larger than 2 or equal to 2. Then it follows that
\footnote{For simplicity, we employ the notation in which
$F(L_D)$ is replaced by $F(D)$.
$D$ represents an
elementary sub-diagram of $L$.}
\begin{equation}
F\Unfoldingsssb=\frac1{e(1)-e(-1)}\left(
e(-\tfrac1{2})F\Pcrossb-e(\tfrac1{2})F\Ncrossb\right)~~~,
\label{suunfoldspinor}
\end{equation}
for the second term on the right hand side of (\ref{tanglecross2})
and
\begin{equation}
F\Unfoldb=\frac1{e(1)-e(-1)}
\left(e(\tfrac1{2})F\Ncrossb-e(-\tfrac1{2})F\Pcrossb\right)~~~,
\label{stunfoldspinor}
\end{equation}
for the second term on the right hand side of (\ref{tanglecross1}).
(\ref{suunfoldspinor}) and (\ref{stunfoldspinor})
are satisfied by tangle operators of oriented links.
\end{lem}
{\it Proof.~}
Let's introduce a colored coupon $D$ \cite{KM} where $D$ is an
isomorphism ${V^2}^*\longrightarrow V^2$ under the assumption $k\geq 2$.
${V^2}^*$ is dual to $V^2$.
Let $G_c$ represent a colored framed graph formed by inserting
the colored coupon $D$ and its inverse $D^{-1}$ on $K$ so that
$p$ extreme points of the link diagram of $L$ are separated.
$K$ represents a link component of $L$
with $k$-color where $k$ represents a $k$-dimensional
representation of $U_q(sl(2,\bf C))$.
There is a fact that
\begin{equation}
F(G_c)=(-1)^{(1-k)p}F(L)~~~.
\label{Dproperty}
\end{equation}
First,  we prove (\ref{stunfoldspinor}).
{}From the property of the quantum $R$ matrix,
the tangle operators of links must satisfy the skein relation:
\begin{equation}
F\Unfoldingb=\frac1{e(1)-e(-1)}
\left(e(\tfrac1{2})F\Pcrossb-e(-\tfrac1{2})F\Ncrossb\right)~~~.
\label{tangleskein}
\end{equation}
When rotated by 90 degrees, it looks like
\begin{equation}
F\Unfoldb=\frac1{e(1)-e(-1)}
\left(e(\tfrac1{2})F\Pcrosssb-e(-\tfrac1{2})F\Ncrosssb\right)~~~.
\label{tangleskein2}
\end{equation}
Moreover, one can rewrite (\ref{tangleskein2}) into the following
form:
\begin{equation}
F\Unfoldb=\frac1{e(1)-e(-1)}
\left(e(\tfrac1{2})F\Ncrossb-e(-\tfrac1{2})F\Pcrossb\right)~~~.
\label{tangleskein3}
\end{equation}
This follows from the following relations:
{\allowdisplaybreaks
\begin{align}
F\Pcrosssb&=F\DDpcrossb=F\Ncrossb~~~,\\
F\Ncrosssb&=F\DDncrossb=F\Pcrossb~~~.
\label{DDequation2}
\end{align}}
Thus we proved (\ref{stunfoldspinor}).
\par
Second, we prove (\ref{suunfoldspinor}).
{}From (\ref{Dproperty}) , it follows that
\begin{align}
&F\Unfoldingsssb=-F\Dunfoldingb\notag\\&=
\frac{-1}{e(1)-e(-1)}\left(e(\tfrac1{2})F\Dncrossb-e(-\tfrac1{2})F\Dpcrossb\right)\\
&=
\frac1{e(1)-e(-1)}\left(
e(-\tfrac1{2})F\Pcrossb-e(\tfrac1{2})F\Ncrossb
\right)~~~.
\notag
\label{Dunfolding2}
\end{align}
In the second equality, we used the skein relation (\ref{tangleskein2}).
In the last equality, we used the following relations:
\begin{equation}
F\Dncrossb=F\Ncrossb,~~~~~
F\Dpcrossb=F\Pcrossb~~~,
\label{DDequation1}
\end{equation}
(\ref{DDequation1}) is immediately derived from (\ref{Dproperty})
because there exist even number of extrema between the two
coupons
$D$ and $D^{-1}$ in both equations of (\ref{DDequation1}).
Thus we completed the proof of the lemma.
\hfill\qed \vspace{0.2cm}\break
We are ready to prove the theorem 4.1.
\hfill\break
{\it Proof~of~Theorem~4.1~}
It is obvious from the skein relation
(\ref{tangleskein}) and the lemma 4.1.
Substitute (\ref{suunfoldspinor}) and (\ref{stunfoldspinor})
into the right hand sides of
(\ref{tanglecross2}) and  (\ref{tanglecross1}) respectively and use
(\ref{tangleskein}).
We arrive at (\ref{resolution}) for the $j=1$ case .
\hfill\qed \vspace{0.2cm}
\break Let's generalize the theorem 4.1 to cases of $j\geq1$.
\begin{thm}
Let $F\left(L^{(j,0)}\right)$ for $j\geq1$ be
the tangle operator of a graph with j rigid vertices given by
(\ref{resolution}).
It satisfies the spinor identity
\begin{equation}
F\left(L^{(j,0)}\scross\right)
=F\left(L^{(j-1,0)}\sunfold\right)-F\left(L^{(j-1,0)}\suunfold\right)~~~~,
\label{tanglejcross2}
\end{equation}
for the case 1 (see the left diagram of Fig.\ref{Spgraph3}), and
\begin{equation}
F\left(L^{(j,0)}\scross\right)
=F\left(L^{(j-1,0)}\sunfold\right)+F\left(L^{(j-1,0)}\stunfold\right)~~~,
\label{tanglejcross1}
\end{equation}
for the case 2 (see the right diagram of Fig.\ref{Spgraph3}).
\end{thm}
{\it Proof.~}
When $j=1$, the theorem was already proved.
So it is enough to assume $j\geq2$.
Suppose that all the rigid vertices of $L^{(j,0)}$ except the j-th one are
resolved.
Applying the definition (\ref{resolution}) to the resolution,
we find on the left hand side of (\ref{tanglejcross2}) or
(\ref{tanglejcross1}),
\begin{equation}
F\left(L^{(j,0)}\scross\right)
=\frac{1}{{\left(e(\tfrac{1}{2})+e(-
\tfrac{1}{2})\right)}^{j-1}}\times
\sum_{s=1}^{2^{j-1}}F\left(L^{(1,0)}_s\scross\right)~~~,
\label{jresolution}
\end{equation}
where $\{L^{(1,0)}_s\scross\}$ is a set of all
graphs with only one rigid vertex obtained by resolving $j-1$ rigid vertices.
It is necessary to distinguish graph types appearing in
(\ref{tanglejcross2}) from those appearing in (\ref{tanglejcross1}) .
Thus $F\left(L^{(j,0)}\scross\right)$ is expressed by a sum of tangle
operators
of $2^{j-1}$ graphs with only one rigid vertex.
\par
On the other hand, we can resolve all $j-1$
rigid vertices
on the right hand sides of (\ref{tanglejcross2}) and
(\ref{tanglejcross1}) in the same way.
After all, they take the following forms:
{\allowdisplaybreaks
\begin{align}
&F\left(L^{(j-1,0)}\sunfold\right)
-F\left(L^{(j-1,0)}\suunfold\right)\notag\\
&=\frac{1}{{\left(e(\tfrac{1}{2})+e(-\tfrac{1}{2})\right)}^{j-1}}\times
\sum_{s=1}^{2^{j-1}}F\left(L^{(0,0)}_s\sunfold\right)
 -F\left(L^{(0,0)}_s\suunfold\right)~~~,
\label{jtanglejcross1}
\end{align}}
for the case 1, and
{\allowdisplaybreaks
\begin{align}
&F\left(L^{(j-1,0)}\sunfold\right)+F\left(L^{(j-
1,0)}\stunfold\right)\notag\\
&=\frac{1}{{\left(e(\tfrac{1}{2})+e(-\tfrac{1}{2})\right)}^{j-1}}\times
\sum_{s=1}^{2^{j-1}}
F\left(L^{(0,0)}_s\sunfold\right)+F\left(L^{(0,0)}_s\stunfold\right)~~~,
\label{jtanglejcross2}
\end{align}}
for the case 2.
We used the definition (\ref{resolution}).
What we only have to do is to show that for every $s$
\begin{equation}
F\left(L^{(1,0)}_s\scross\right)=F\left(L^{(0,0)}_s\sunfold\right)-
F\left(L^{(0,0)}_s\suunfold\right)~~~,
\label{jstanglejcross1}
\end{equation}
for the case 1, and
\begin{equation}
F\left(L^{(1,0)}_s\scross\right)=
F\left(L^{(0,0)}_s\sunfold\right)+F\left(L^{(0,0)}_s\stunfold\right)~
{}~~,
\label{jstanglejcross2}
\end{equation}
for the case 2.
The graph $L^{(1,0)}_s\scross$ ($L^{(1,0)}_s\scross$)
and the links $L^{(0,0)}_s\sunfold$ and $L^{(0,0)}_s\suunfold$
($L^{(0,0)}_s\sunfold$ and $L^{(0,0)}_s\stunfold$) in
(\ref{jstanglejcross1}) ((\ref{jstanglejcross2})) are identical
outside the subdiagrams depicted
for each $s$, if all the orientations of edges and loops are neglected.
Thus the proof of (\ref{jstanglejcross1}) and (\ref{jstanglejcross2})
is accomplished by repeating the proof of the theorem 4.1.
\hfill\qed \vspace{0.2cm}
\break
\par
The inverse problem is also of interest to us.
It is whether
operators $F(L^{(j,0)})$ characterized
by (\ref{tanglejcross2}) and (\ref{tanglejcross1}) must be subject to
(\ref{resolution}) or not.
We can easily resolve it.
As a consequence, it follows that the
operators subject to the spinor identity must be subject to (\ref{resolution}).
\par
Thus we constructed the tangle operators $F(L^{(k,j-k)})$
which are uniquely identified with the CS vacuum expectation values of Wilson
loops with  j transverse double points.
The Casimir operators are inserted at j-k transverse double
points one by one.
Of course, equivalence of  $F(L^{(k,j-k)})$ and $Z(L^{(k,j-k)})$
should be verified order by order with respect
to the CS coupling constant $k$ in the perturbation theory.
There is a technical limit to achieve it.
Computation of terms in higher orders is extremely difficult.
This is one of motivations of the present paper.
\par
The following is now our main conclusion.
\begin{cor}
Let $P\Lonecroossb$ be the graph invariant of Vassiliev type
given by a graph with only one rigid vertex.
Then it can be expressed by the two types of tangle operators of graphs:
{\allowdisplaybreaks
\begin{align}
P\Lonecroossb&\equiv P\Lzpcrossb-P\Lzncrossb \notag\\
&=\alpha^{-\omega(L^{(1)})}(e(1)+e(-1)-1)(e(1)-e(-1))\notag\\
                    &\phantom{=}
                        \times[2F\left(L^{(0,1)}\sccross\right)
 -\frac{1}{2}
               \frac{e(1)+e(-1)+1}{e(1)+e(-1)-1}
               F\left(L^{(1,0)}\scross\right)]~~.
\label{cor2}
\end{align}}
Let's call the expression of
(\ref{cor2}) {\bf the Casimir insertion representation of the
graph invariants of Vassiliev type} in the context of the Reshetikhin-Turaev
construction of quantum group topological invariants.
\end{cor}
{\it Proof.~}
Solving (\ref{resolution})
and (\ref{tanglecasimir}) for $j=1$,we immediately find
{\allowdisplaybreaks
\begin{align}
F\left(L^{(0,0)}\spcross\right)
&=\tfrac{1}{2}\left(e(\tfrac{1}{2})
                 +e(-\tfrac{1}{2})\right)F\left(L^{(1,0)}\scross\right)
                  +2\left(e(\tfrac{1}{2})-e(-\tfrac{1}{2})\right)
                       F\left(L^{(0,1)}\sccross\right)
{}~~~,\notag\\F\left(L^{(0,0)}\sncross\right)&=\tfrac{1}{2}\left(e(\tfrac{1}{2})
+e(-\tfrac{1}{2})\right)F\left(L^{(1,0)}\scross\right)
                  -2\left(e(\tfrac{1}{2})
                    -e(-\tfrac{1}{2})\right)
                        F\left(L^{(0,1)}\sccross\right)~~~.
\label{sol1}
\end{align}}
According to
\begin{equation}
P\Lzpcrossb=
\alpha^{-\omega\left(L^{(0,0)}\spcross\right)}
F\left(L^{(0,0)}\spcross\right)~, \quad
P\Lzncrossb=
\alpha^{-\omega\left(L^{(0,0)}\sncross\right)}
F\left(L^{(0,0)}\sncross\right)~~~,
\label{sol2}
\end{equation}
we arrive at (\ref{cor2})
after substitution of (\ref{sol1}) into (\ref{cor2}).
\hfill\qed \vspace{0.2cm}
\break
The restriction to $j=1$ is not crucial.
More general cases of $j>1$ can be considered in the same way
and will be considered soon later.
But then, expressions are more complicated.
\par
We shall end this section with a few conclusions derived from the
corollary 4.1.
First, from the Casimir insertion representation (\ref{cor2}) ,
it is obvious that the tangle operator of a rigid vertex is nothing but
\begin{equation}
\Fvertex\equiv\left(e(1)+e(-1)-1\right)\left(e(1)-e(-1)\right)
\left(2\Fccross-\frac1{2}\frac{e(1)+e(-1)+1}{e(1)+e(-1)-
1}\Fcross\right)~~~.
\label{Vertex}
\end{equation}
where
{\allowdisplaybreaks
\begin{align}
\Fcross&=\frac{1}{e\left(\tfrac1{2}\right)+e\left(-\tfrac1{2}\right)}
\left(\Fpcross+\Fncross\right)~~~,\\
\Fccross&=
\frac{1}{4\left(e\left(\tfrac1{2}\right)-e\left(-\tfrac1{2}\right)\right)}
\left(\Fpcross-\Fncross\right)~~~.
\end{align}}
$\Fpcross$ ($\Fncross$) is
the tangle operator of the positive (negative) crossing.
(\ref{Vertex}) says that $\Fvertex$ is a linear combination of tangle operators
of
two kinds of rigid vertices, i.e., $\Fccross$ and $\Fcross$.
\par
Second, the graph invariant of Vassiliev type $P\Lonecroossb$
can immediately be generalized to graph invariants $P\Ljcroossb$ for $j>1$.
Then it is obvious that
\begin{equation}
P\Ljcroossb=\alpha^{-\omega(L^{(j)})}
{\tilde F}(\bar L)\circ\otimes^j_{t=1}\Ftvertex~~~,
\label{jvertex}
\end{equation}
where ${\tilde F}(\bar L)$ represents a tangle operator given by a
complement $\bar L$ obtained by cutting neighborhoods of the
j rigid vertices out of the graph $L^{(j)}$.
$\Ftvertex$ represents a tangle operator of the t-th rigid
vertex.
\par
Third, it is interesting to compare our formula (\ref{cor2}) with Kauffman's
approximate one computed in \cite{KAUFF3} .
He found the Casimir insertion representation
of the graph invariants of Vassiliev type based on
the perturbative analysis of the CS quantum gauge field theory.
He obtained  the following approximate formula
\begin{equation}
P\Lonecroossb=\alpha^{-\omega(L^{(j)})}\left(\tfrac{4\pi
i}{k}\right)\left(Z\left(L^{(1,0)}\scross\right)-\tfrac{3}{4}Z\left(L^{(0,1)}\sccross\right)
\right)
\end{equation}
in the leading order of the CS coupling constant $k$.
It is easy to check that this coincides with our formula (\ref{cor2}).
We can think of the expression ($\ref{cor2}$) as a non-perturbative expression
of
Kauffman's approximate formula.
We must mention a fact that
Kauffman's approximate formulae can exist associated with any
Lie algebra.
A generalization of our argument of the present paper restricted to
$U_q(sl(2,\bf C))$
to any quasi-triangular Hopf algebra remains to be investigated.
\subsection{Locally Integrable Condition and 6-Valent Graph Invariants}%%%%%%%%
In the theory of the Vassiliev invariants,
the Vassiliev invariants must be subject to the local integrability
condition (which is often called the four term relation)
for resolution of a transverse triple point.
The transverse triple point is allowed to form when there are more than
one double point.
The local integrability condition has been discussed to a great extent by
Birman and
Lin\cite{BL1}\cite{Lin1}.
It is the following relation among the Vassiliev invariants of four singular
links
obtained by moving one of three axes
along the other two axes in the vicinity of the transverse
triple point:
\begin{equation}
v_i\left(\Northb\right)-v_i\left(\Southb\right)
+v_i\left(\Eastb\right)-v_i\left(\Westb\right)=0~~~,
\label{fourterm}
\end{equation}
where $v_i$ represents any one of the Vassiliev invariants of order $i$.
The Vassiliev invariants of order $i$ appear
as coefficients in an expansion of the Jones polynomial of a link with
respect to $x$ where $q=e^x$.
In accordance with (\ref{fourterm}), we can write
the local integrability condition satisfied by the graph
invariants of Vassiliev type defined in \S2.
It is
\begin{equation}
P\left(\North\right)-P\left(\South\right)
+P\left(\East\right)-P\left(\West\right)=0~~~.
\label{fourtermgraph}
\end{equation}
Each term on the left hand side of (\ref{fourtermgraph})
is given by the tangle operator of a graph.
(\ref{fourtermgraph}) implies existence of a graph invariant
given by the 6-valent vertex graph.
We can define it as follows.
\begin{defn}
Suppose that
we are given graph invariants of Vassiliev type given by graphs
with 4-valent
vertices.
We define the graph invariant of Vassiliev type given by the 6-valent vertex
graph
by
{\allowdisplaybreaks
\begin{align}
P\left(\underbrace{\triple}_{G}\right)&=
P\left(\underbrace{\North}_{G_1}\right)
-P\left(\underbrace{\South}_{G_2}\right)\notag\\
&=P\left(\underbrace{\West}_{G_3}\right)
-P\left(\underbrace{\East}_{G_4}\right)~~~.
\label{Triple}
\end{align}}
\end{defn}
The graph invariant given by the 6-valent vertex  is a gap
between the graph invariant of $G_1$ ($G_3$) and the graph invariant of $G_2$
($G_4$).\par
In this sub-section, we attempt to show
how (\ref{fourtermgraph}) is satisfied in the Casimir insertion representation
in which
the graph invariants of Vassiliev type are identified with
the vacuum expectation values of Wilson loops with transverse intersection
points.
We proceed assuming such identification in the followings.
Let's use the Kauffman bracket $Z(G)$ instead of the tangle operator $F(G)$.
In the work \cite{KAUFF3},
Kauffman noticed the following fact based on the perturbative analysis of the
CS quantum gauge field theory.
It said that
(\ref{fourtermgraph}) can be satisfied by the
commutation relations of a Lie
algebra\cite{KAUFF3} as far as only the Vassiliev invariants in the
top raw of the actuality table\cite{Bir1}\cite{BL1} are concerned.
But we wonder if such a point of view can work in general, i.e., even for
all the Vassiliev invariants.
We attempt to verify that the local integrability condition is satisfied by all
the
Vassiliev invariants using the Casimir insertion
representation.
\par
According to the definition (\ref{Triple}), the Casimir insertion
representation
of graph invariants of the 6-valent vertex graphs
is defined as follows.
The 6-valent vertex can be regarded as the transverse triple point formed
by two double points.
{}From the first equality of (\ref{Triple}), it is given by
{\allowdisplaybreaks
\begin{align}
&P\left(\underbrace{\triple}_{G}\right)=\alpha^{-\omega(G_1)}
Z\left(\underbrace{\nvvb}_{G_1}\right)
-\alpha^{-
\omega(G_2)}Z\left(\underbrace{\svvb}_{G_2}\right)\notag\\
&\phantom{P\left(\underbrace{\triple}_{G}\right)}
=\alpha^{-\omega(G_1)}Z\left(\nvv\right)-\alpha^{-\omega(G_2)}Z\left(\svv\right)~~~.
\label{sixgraphcasim1}
\end{align}}
We formed a disk spanned by three axes applying (\ref{sol1}).  We can use a
fact that
arbitrary shift of any one of axes forming the disk within the disk doesn't
make any change of value of the graph invariant.
It is natural
from a point of view of the CS path integral.
Thus we define the Casimir insertion representation
of the 6-valent graph by the
following limit:
{\allowdisplaybreaks
\begin{align}
&P\left(\triple\right)
=lim_{\text{Disk}\rightarrow\text{point}}\notag\\
&
a\alpha^{-
\omega(G)}Z\left(\nvvs\right)
+b\alpha^{-
\omega(G)}Z\left(\nvvc\right)~~~,\notag\\
&
-c\alpha^{-\omega(G)}
Z\left(\svvs\right)
-d\alpha^{-
\omega(G)}Z\left(\svvc\right)~~~.
\label{sixgraphcasimir1}
\end{align}}
where
{\allowdisplaybreaks
\begin{align}
a&=\frac1{2}\left(e(-1)+e(-2)\right),\quad
b=2\left(e(-1)-e(-2)\right)~~~,\notag\\
c&=\frac1{2}\left(e(1)+e(2)\right),
\quad  d=-2\left(e(2)-e(1)\right)~~~.
\end{align}}The disk formed by three axes shrinks to a point in the limit.
On the other hand, from  the second equality of (\ref{Triple}),
the Casimir insertion representation of the 6-valent graph is given by
{\allowdisplaybreaks
\begin{align}
P\left(\triple\right)
&=\alpha^{-\omega(G_3)}Z\left(\wvvb\right)
-\alpha^{-\omega(G_4)}Z\left(\evvb\right)
\\
&
=\alpha^{-\omega(G_3)}Z\left(\wvv\right)
-\alpha^{-\omega(G_4)}Z\left(\evv\right)
\notag\\
&=lim_{\text{Disk}\rightarrow\text{point}}\notag\\
&~~~~~~~~~~~~~~~~a\alpha^{-\omega(G)}Z\left(\wvvs\right)
+b\alpha^{-
\omega(G)}Z\left(\wvvc\right)~~~,\notag\\
&~~~~~~~~~~~~~~~~~~
-c\alpha^{-\omega(G)}
Z\left(\evvs\right)
-d\alpha^{-
\omega(G)}Z\left(\evvc\right)~~~.
\label{sixgraphcasimir2}
\end{align}}
The deformations of the vertical axis to form disks make it easy to
compare (\ref{sixgraphcasimir1} ) and (\ref{sixgraphcasimir2}).
We can easily check  that they coincide.
The following is available in understanding how to compute them explicitly.
\par
As we have seen above, the Casimir insertion representation of the 6-valent
graph
invariant is determined by the limit of shrinking the disk to a point.
Let's calculate it more explicitly as an exercise.
Recall that the Casimir insertion representation of
the tangle operator of the rigid vertex is given by (\ref{Vertex}):
\begin{equation}
\Fvertex=C_1\Fcross+C_2\Fccross~~~,
\label{Vertex1}
\end{equation}
where
{\allowdisplaybreaks
\begin{align}
C_1&=a-c=-\frac1{2}\left(e(1)-e(-1)\right)\left(e(1)+e(-1)+1\right),\notag\\
C_2&=b-d=2\left(e(1)-e(-1)\right)\left(e(1)+e(-1)-1\right)~~~.
\end{align}}
The rigid vertex is replaced by a transverse double point of the Wilson loops
at the moment.
We can apply it to the computation.
After some algebra, we find
%%%%%%%
{\allowdisplaybreaks
\begin{align}
P\left(\triple\right)&=
P\left(\North\right)-P\left(\South\right)\notag\\
&=aZ\left(\nvvsb\right)+bZ\left(\nvvcb\right)\notag\\
&\phantom{=aP\left(\nvvsb\right)}
-cZ\left(\svvsb\right)-dZ\left(\svvcb\right)\notag\\
&=H_1+H_2+H_3+H_4+H_5~~~.
\label{TripleNS}
\end{align}}
where
{\allowdisplaybreaks
\begin{align}
&H_1
=(C_1)^3\alpha^{-\omega(G)}Z\left(\trib\right)~~~,\\
&H_2=(C_1)^2C_2\\
&\times\left\{
\alpha^{-\omega(G)}Z\left(\triplea\right)
+\alpha^{-\omega(G)}Z\left(\tripleb\right)
+\alpha^{-\omega(G)}Z\left(\tripled\right)
\right\}~~~,\notag\\
&H_3=C_1(C_2)^2\\
&\times\left\{
D^{ab}_{cd}\alpha^{-\omega(G)}Z\left(\triplen\right)
+D^{ab}_{cd}\alpha^{-\omega(G)}Z\left(\triplee\right)
+D^{ab}_{cd}\alpha^{-\omega(G)}Z\left(\triplef\right)
\right\}~~~,\notag\\
&H_4=\left(C_1C_2(b+d)+\frac1{2}(C_2)^2(a+c)\right)
i\epsilon_{abc}\alpha^{-\omega(G)}Z\left(\tripleabc\right)~~~,\\
&H_5=
(C_2)^2\left\{
b\alpha^{-\omega(G)}Z\left(\tripleg\right)
-d\alpha^{-\omega(G)}Z\left(\tripleh\right)\right\}~~~.
\end{align}}
where we used the following symmetrizer
$D^{ab}_{cd}=\tfrac1{2}\left(\delta^a_c\delta^b_d+\delta^a_d\delta^b_
c\right)$.
$\epsilon_{abc}$ is the structure constant of the Lie algebra, i.e.,
$[T_a,T_b]=i\epsilon_{abc}T_c$.
\par
Let's end this sub-section with the following result.
\begin{rem}
We found the Casimir insertion representation of the Vassiliev
invariants of graphs with 4-valent and 6-valent vertices.
The local integrability condition
(or the four term relation) can be satisfied by the commutation relation of the
Lie
algebra.
\end{rem}
The point in the remark has been pointed out by Bar-Natan,
Kauffman et. al. from a point of view of
the perturbation theory of the CS quantum gauge field
theory.
Thus the present work can provide a non-perturbative perspective on
their observation.
%
%
%
%
%%%%%%%%%%%%%%%%%%%%%%%%%%%%%%%%%%%%%%%%%%%%%%%%%
%%%%%%%%%   Section 5      %%%%%%%%%%%%%%%%%%%%%%
\section{Application to Physics}%%%%%%%%%%%%%%%%%
%%%%%%%%%%%%%%%%%%%%%%%%%%%%%%%%%%%%%%%%%%%%%%%%%
%%%%%%%%%%%%%%%%%%%%%%%%%%%%%%%%%%%%%%%%%%%%%%%%%
\par
In the previous sections, we discussed quantum group invariants
of rigid vertex graphs with four- and six-valent rigid vertices and found
 the Casimir insertion representation of them.
It is an analytic expression in the CS quantum gauge field theory
for the graph invariants. Importance of the Casimir insertion representation
comes from mathematical and physical points of view.
{}From the mathematical point view, it fills out a gap between an algebraic
definition of the graph invariants based on the quantum groups and an analytic
definition of them based on the path integral
quantization of the CS field theory.
It shed light on the local integrability conditions to
which the graph invariants of Vassiliev type must to be subject.
On the other hand, from the physical point of view,
it enables us to discuss whether the Vassiliev invariants can be
physical states of the 4D quantum gravity of Ashtekar or not.
The Casimir insertion representation plays a crucial role in such a discussion.
\par
In this section, we are interested in the latter.
Ashtekar's gravity provides a scheme of non-perturbative analysis in
the canonical quantization.
The loop space representation of wave functions in Ashtekar's gravity
is most significant ,in which wave functions
are defined over a loop space, i.e., a space of maps from $S^1$ to a 3-space.
We here aim at clarifying that wave functions given by the graph invariants of
Vassiliev type can be physical states in the loop space representation of the
quantum gravity.
Let's begin with a brief review on the canonical quantization of Ashtekar's
gravity and the loop space representation of wave functions.
\subsection{A Brief Review on the Loop Space Representation in the Quantum
Gravity}%%%%%%%%
\par
Let $M^4=\Sigma^3\times \bf R$
be a (real analytic) 4-dimensional differential manifold with a co-dimension
one
foliation,
and $\Sigma^3(t)$ a leaf. $t$ is a parameter of time. It is given by
$t=\tau(\Sigma^3)$ in
terms of a smooth map $\tau:~\Sigma^3\rightarrow\bf R$.
We suppose that we are given complex-valued functionals $\psi(A:\Sigma^3(t))$
defined
over an affine space $\cal A$
\footnote{
In Ashtekar's gravity,
the self-dual connection $A^4$ ($A^4=A_0dt+A_idx^i$) and the tetrad
$\tilde E$ defined over $M^4$ are dynamical variables.
In the (3+1)-decomposition and in $A_0=0$ gauge, $A\equiv A_idx^i$ is
a coordinate of the configuration space on which
wave functions are defined in the canonical quantization.
}
 of $su(2)$-valued connection 1-forms over $\Sigma$
\footnote{For brevity, we replace $\Sigma^3(t)$ by $\Sigma$.}
{}.
They are sections of a line bundle over $\cal A$ specified by a set of
constraints.
In the canonical quantization, the constraints of Ashtekar's quantum gravity
with
vanishing cosmological constant take the following forms:
 {\allowdisplaybreaks
\begin{align}
\hat{\frak G}[\epsilon^i]\psi(A:\Sigma)
&=\int_{\Sigma}d^3x\epsilon^i(x){\hat {\cal G}}_i\psi(A:\Sigma)\notag\\
&=i\int_{\Sigma}d^3x
\epsilon^i(x)\cal D_a\frac{\delta}{\delta {A^i_a(x)}}\psi(A:\Sigma)=0~~~,
\label{const1}\\
\hat{\frak M}[N^a]\psi(A:\Sigma)&=\int_{\Sigma}d^3xN_a (x){\hat {\cal M}}^a
\psi(A:\Sigma)\notag\\
&=i\int_{\Sigma}d^3xN_b (x)
\frac{\delta}{\delta A^i_a(x)}F^i_{ab}\psi(A:\Sigma)=0~~~,
\label{const2}\\
\hat{\frak H}[N]\psi(A:\Sigma)&=\int_{\Sigma}d^3x\underset{\sim}{N}(x)
\hat{\cal H}\psi(A:\Sigma)\notag\\
&=\int_{\Sigma}d^3x \underset{\sim}{N}(x) \epsilon^{ijk}\frac{\delta}{\delta
A^i_a(x)}
\frac{\delta}{\delta
A^j_b(x)}F^k_{ab}\psi(A:\Sigma)=0~~~,
\label{const3}
\end{align}
}
where $\epsilon^i(x)$, the shift functions $N_a(x)$ and the lapse function
$\underset{\sim}{N}(x)$
are analytic functions over $M^4$.
The first constraint is called the Gauss law constraint, the second the
momentum
constraint which generates diffeomorphisms of $\Sigma$,
the last one the Hamiltonian constraint
(which is often called
the Wheeler-DeWitt equation) which generates diffeomorphisms in the time
direction.
\par
Let's introduce the loop space representation of the wave functions.
Remember that the Wilson loops $\cal W(\cal A:G)$ introduced in \S2 are
gauge invariant
objects.
$G$ represents spatial graphs identified with singular links only with
transverse intersection points such as transverse double and triple points.
We introduce wave functions denoted by $\psi(G:\Sigma)$.
They are defined by functional integration
$\psi(G:\Sigma)=\int_{\cal A}d\mu(\cal A)\cal W(\cal A:G)\psi(A:\Sigma)$.
Constraints in the loop space representation
are induced by the constraints (\ref{const1}), (\ref{const2}) and
(\ref{const3})
via partial integration.
To be concrete,
for differential operators $\hat{\cal O}=\hat{\cal M}$ or $\hat{\cal H}$
(their adjoint operators $\hat{\cal O}^\dagger$),
the constraints in the loop space representation are given by
{\allowdisplaybreaks
\begin{align}
\hat{\cal O}_L\psi(G:\Sigma)&=\int_{\cal A}d\mu(\cal A)W(\cal A:G)
\left(\hat{\cal O}\psi(A:\Sigma)\right)\notag\\
&=\int_{\cal A}d\mu(\cal A)\left(\hat{\cal O}^\dagger\cal W(\cal
A:G)\right)\psi(A:\Sigma)=0~~~.
\label{loopconstraint}
\end{align}}
We used partial integration to obtain the second equality.
The adjoint operators are dependent on the measure.
Lack of the Gauss law constraint comes from  a fact that
 the loops space representation is defined to be manifestly gauge invariant by
definition.
\par
It is known that
there is an important class of solutions to
the constraint (\ref{loopconstraint}) in the loop space
representation of the quantum gravity with non-zero cosmological constant.
Let's consider wave functions $\psi_{cs}(A:\Sigma)$ given by the
Chern-Simons action $S_{cs}(A)$, i.e.,
$\psi_{cs}(A:\Sigma)\equiv\exp(-\frac{6}{\Lambda}\int_{\Sigma} tr\left(
A\wedge dA+\frac{2}{3}A\wedge A\wedge A
\right))
=\exp(ikS_{cs}(A))$
where $S_{cs}(A)=\frac{1}{4\pi}\int_{\Sigma} tr\left(A\wedge
dA+\frac{2}{3}A\wedge A\wedge A
\right)$ and $\Lambda=i\frac{24\pi}{k}$ (the cosmological constant).
We call them CS states.
This type of solutions was initially discussed by Kodama\cite{Kodama1}.
The wave functions in the loops space representation are just Kauffman's
brackets.
One can easily check that they satisfy the Hamiltonian constraint with non-zero
cosmological constant :
${\hat {\cal H}}^\Lambda_L\psi_{cs}(G:\Sigma)=\left({\hat {\cal H}}_L
-\frac{\Lambda}{6}\hat{det}(\tilde E)_L\right)\psi_{cs}(G:\Sigma)=0$.
${\hat {\cal H}}_L$ is induced by the Hamiltonian operator with non-zero
cosmological constant:
\begin{equation}
{\hat {\cal H}}^\Lambda\equiv\epsilon^{ijk}
\frac{\delta}{\delta A^i_a}\frac{\delta}{\delta A^j_b}F^k_{ab}
+\frac{\Lambda}{6}\epsilon_{abc}\epsilon^{ijk}
\frac{\delta}{\delta A^i_a}
\frac{\delta}{\delta A^j_b}
\frac{\delta}{\delta A^k_c}~~~.
\label{Hamiltonian}
\end{equation}
The reason why the CS states satisfy the Hamiltonian constraint is obvious from
a fact:
\begin{equation}
\frac{\delta}{\delta A^i_a}\psi_{cs}(A:\Sigma)
=i\frac{k}{8\pi}\epsilon_{abc}F^i_{bc}\psi_{cs}(A:\Sigma)~~~.
\end{equation}
In showing ${\hat {\cal H}}^\Lambda_L\psi_{cs}(G:\Sigma)=0$,
we assumed that partial integration is possible in the CS functional
integral.
%%%%%%%%%%%%%%%%%%%%%%%%%%%%%%%%%%%%%%%%%%%%%%%%%%%%%%%%%%%%%%%%%%%%%%
%%%%%%%%%%%%%%%%%%%%%%%%%%%%%%%%%%%%%%%%%%%%%%%%%%%%%%%%%%%%%%%%%%%%%%
\subsection{Spin-Network States and the Hamiltonian Constraint}%%%%%%%
%%%%%%%%%%%%%%%%%%%%%%%%%%%%%%%%%%%%%%%%%%%%%%%%%%%%%%%%%%%%%%%%%%%%%%
%%%%%%%%%%%%%%%%%%%%%%%%%%%%%%%%%%%%%%%%%%%%%%%%%%%%%%%%%%%%%%%%%%%%%%
\par
The q-analog of Penrose's spin-network was initially discussed by
Kauffman\cite{KAUFF4}.
He arrived at the quantum group graph invariants of Dubrovnik type
(a graph extension of the Kauffman polynomials of links), that is,
the quantum group invariants of unoriented graphs.
However,
in \S 2, we saw that the spinor identity which is indispensable to the q-analog
of
Penrose's spin-network requires the orientation of graphs.
In \S 3, following the Reshetikhin-Turaev construction of 3-manifold
invariants,
we defined the q-analog of Penrose's spin-network taking into account of the
orientation
in an algebraic manner. The spinor identity plays a role of resulting vertices
of graphs.
The discussion of the graph invariants in the context of
the q-analog of Penrose's spin-network led us to find the analytic expression
(i.e., the Casimir insertion representation) for the graph invariants of
Vassiliev type.
Thus we arrived at a graph extension of the Jones polynomials.
It is likely that a graph extension of the HOMFLY polynomials will be
obtained by discussing Mandelstam's identity for $SU(N)$.
\par
We here show that
the graph invariants of Vassiliev type play a crucial role in the canonical
quantum gravity of
Ashtekar\footnote{Importance of the Vassiliev invariants
in Ashtekar's quantum gravity was first emphasized by J.Baez\cite{Baez1}.}.
They are regarded as the physical states of the universe.
Let's call them spin-network states.
To be precise,
they satisfy the Hamiltonian constraint with vanishing cosmological constant.
The following is a conclusion of the present paper.
\begin{thm}
Let $P(G:\Sigma)$ be a graph invariant of Vassiliev type.
G represents a singular link only with transverse double and triple points.
Then it satisfies
\begin{equation}
{\hat{\cal H}}_LP(G:\Sigma)=0~~~.
\end{equation}
\end{thm}
{\it Proof.~}
It suffices to prove ${\hat{\cal H}}_LP(L:\Sigma)=0$ because $P(G:\Sigma)$ can
be obtained
in terms of link invariants $P(L:\Sigma)$ by definition.
It is obvious from the following.
Let's expand the link polynomials $P(L:\Sigma)$
 with respect to the inverse of the CS coupling constant $\frac1{k}$
or the cosmological constant $\Lambda$.
Let $f_n(L)$ be a coefficient in the n-th order.
It is given by a product of a group factor and a geometric
factor (or an analytic factor)\cite{AL}.
The Hamiltonian operator acts on the latter factor.
As far as the link is concerned,
the action of the Hamiltonian operator ${\hat{\cal H}}_Lf_n(L)$
must vanish for all n because it is able to have supports only at intersection
points or
at kinks\cite{BGP}\cite{Brug}.
The Hamiltonian operator generates shifts of line segments along tangent
vectors
at the intersection points or at the kinks accompanied with recomposition
of loops and partial orientation reverse.
The Hamiltonian operator doesn't generate any shifts when acting
on a functional of loops with no intersection points and no kinks.
Thus it follows that ${\hat{\cal H}}_LP(G:\Sigma)=0$.
\hfill\qed \vspace{0.2cm}
\break
Another proof comes from a fact that invariants associated to graphs with
transverse triple points are calculated by the limit of shrinking a area
surrounded by three axes to a point on a 2-dimensional disk.
(We clarified it in  \S 3. )
But the action of the Hamiltonian operator on functionals over the space of
singular links is non-trivial only at transverse triple points.
So the action of the Hamiltonian operator inevitably vanishes.
%%%%%%%%%%%%%%%%%%%%%%%%%%%%%%%%%%%%%%%%%%%%%%%%%%%%%%%%%%%%%%%
%%%%%%%%%%%%%%%%%%%%%%%%%%%%%%%%%%%%%%%%%%%%%%%%%%%%%%%%%%%%%%%
\subsection{Spin-network States and Half Flat Geometry}%%%%%%%%
%%%%%%%%%%%%%%%%%%%%%%%%%%%%%%%%%%%%%%%%%%%%%%%%%%%%%%%%%%%%%%%
%%%%%%%%%%%%%%%%%%%%%%%%%%%%%%%%%%%%%%%%%%%%%%%%%%%%%%%%%%%%%%%
\par
The spin-network states $\psi_{cs}(G:\Sigma)=P(G:\Sigma)$
are physical states in a sense that
they satisfy all the constraints of the quantum gravity of Ashtekar.
Our interest is physical implication of them.
Remember that spin-network states are defined by functional integral of
the Wilson loops
with the Casimir insertion at their intersection points.
In general, in the Chern-Simons path integral, i.e.,
contribution of flat connections over $\Sigma$ dominates in the spin-network
states.
Flatness of the self-dual curvature $F(A)$ implies the half-flatness of the 4-
dimensional geometry in Ashtekar's gravity\cite{Ash2}. In the connection
formalism of
Ashtekar's gravity
with vanishing cosmological constant, Ricci flatness guarantees that
the geometry is hyperk\"{a}hler\cite{Sal1}.
The hyperk\"{a}hler structure is characterized by the Ricci flatness and
three complex structures $I$, $J$ and $K$ with the
quatanionic property, i.e., $I^2=J^2=K^2=-1$ and $IJ=-JI=K$.
\par
Let's consider the complex structures on the loop space representation of the
quantum gravity. First, in a case in which loops have no intersection,
it is well-known that a complex structure can be introduced in the normal
bundle over a space of loops (see Fig.\ref{Frenet}).
%%%%%%figure%%%%%%%%%%%%%%%%%%%%%%
\begin{figure}
\begin{center}\epsfile{file=Frenet,scale=0.4}
\end{center}
\caption{The Fren\(\acute e\)t frame and the complex Structure}
\label{Frenet}
\end{figure}
%%%%%%%%%%%%%%%%%%%%%%%%%%%%%%%%%%%%%%%%%%%%%%%%%%%%%%
Given such a complex structure, the Fren\(\acute e\)t frame along the loops is
fixed.
Let $N$ be a section of the normal bundle along a loop,
then another vector $B$ is defined by $B=T\times N=J_pN$.
$T$ is a section of the tangent bundle.
The Fren\(\acute e\)t frame is given by the three vectors $T$, $N$ and $B$.
Our interest is a case of intersecting links with triple points.
At each transverse triple intersection point, it is likely that we can define
three almost complex structures.
The three almost complex structures are introduced in three normal bundles
given by three independent tangent vectors at the intersection point.
However, they can not be given to be independent because the Fren\(\acute e\)t
frames given by the three
almost complex structures can
coincide by rotations of frames and a redefinition of almost complex
structures.
Without lack of generality,
we can put a situation that
$\tilde B=I_p\tilde T$, $\tilde N=J_p\tilde B$ and $\tilde T=K_p\tilde N$.
One can easily check that
the almost complex structures $I$, $J$ and $K$
must satisfy a set of compatibility conditions which is just a set of
quaternionic relation.
Such a property as the hyperk\"{a}hler structure on the space of singular
links should be further explored.
\par
Thus we obtain the following consequence:
\begin{rem}
In the canonical quantization of Ashtekar's gravity with
vanishing cosmological constant, the spin-network states given by
the graph invariants of Vassiliev type are solutions to the Hamiltonian
constraint.
Furthermore, they are dominated by the half flat geometry in
the classical limit.
\end{rem}
\par
Wave functions that we have discussed so far don't describe the physically
realistic universe.
The reason is that non-degenerate metrics are given only at the transverse
triple
points which are discrete.
An idea of constructing wave functions to describe the realistic universe stems
from
the following inductive limit of the Vassiliev invariants,
$\varinjlim E^{(-j.j)}_\infty$.
Wave functions defined over
an open dense set of transverse triple points must belong to such a group
and are expected to describe the physically realistic universe.
\par
We finally aim at construction of the quantum Hilbert space of
quantum gravity of Ashtekar in our study.
It may be considered either in the connection formalism or in the loop space
representation. It will be constructed by a few steps.
First, we  define wave functions to be tri-holomorphic with respect to
the three complex structures of the hyperk\"{a}hler structure.
Second, we introduce
a physical inner product given by
a hermitian pairing to be compatible with the constraints of the quantum
gravity, i.e.,
the momentum constraint and the Hamiltonian constraint.
The idea of geometric realization of the rational conformal field
theories\cite{Ax1}\cite{Ax2}\cite{Hit1}
will be also available in the quantization of Ashtekar's gravity.
Namely,the Hamiltonian constraint will be closely related to
the Kodaira-Spencer class associated
to deformations of the hyperk\"{a}hler structure.
These things will be discussed elsewhere. %
\vspace{1cm}
\hfill\break
{\it Acknowledgment.~~}\hfill\break
The author wishes to thank J. Baez, H.Kodama, L.Smolin and T.Yoneya for
giving him helpful comments.
This work is supported by Grant-in-Aid for Scientific
Research from  the Ministry of Education, Science and Culture of
Japan.
\end{document}